\documentclass{article}

\usepackage{amssymb,amsmath,amsfonts}
\usepackage{algorithm,algpseudocode}
\usepackage{graphicx,here}
\usepackage{subfigure}
\usepackage{epsfig}
\usepackage{multirow}
\usepackage{booktabs}
\usepackage{url}
\usepackage{cite}
\usepackage{epstopdf}

\begin{document}

\title{Multiobjective optimization to a TB-HIV/AIDS coinfection optimal control problem\thanks{This 
is a preprint of a paper whose final and definite form is with 
'Computational and Applied Mathematics', ISSN 0101-8205 (print), ISSN 1807-0302 (electronic). 
Submitted 04-Feb-2016; revised 11-June-2016 and 02-Sept-2016; accepted for publication 15-March-2017.}}

\author{Roman Denysiuk$^1$\\
\texttt{roman.denysiuk@algoritmi.uminho.pt}
\and
Cristiana J. Silva$^2$\\
\texttt{cjoaosilva@ua.pt}
\and
Delfim F. M. Torres$^2$\\
\texttt{delfim@ua.pt}}


\date{$^1$Algoritmi R\&D Center, University of Minho, Portugal\\[0.3cm]
$^2$Center for Research and Development in Mathematics\\ and Applications (CIDMA),
Department of Mathematics,\\ University of Aveiro, 3810-193 Aveiro, Portugal}

\maketitle


\begin{abstract}
We consider a recent coinfection model for Tuberculosis (TB), 
Human Immunodeficiency Virus (HIV) infection and Acquired Immunodeficiency 
Syndrome (AIDS) proposed in [Discrete Contin. Dyn. Syst. 35 (2015), no.~9, 4639--4663].
We introduce and analyze a multiobjective formulation 
of an optimal control problem, where the two conflicting 
objectives are: minimization of the number of HIV infected 
individuals with AIDS clinical symptoms and coinfected
with AIDS and active TB; and costs related to prevention
and treatment of HIV and/or TB measures. The proposed approach eliminates
some limitations of previous works. The results of the numerical study provide
comprehensive insights about the optimal treatment policies and the population
dynamics resulting from their implementation. Some nonintuitive conclusions
are drawn. Overall, the simulation results demonstrate the usefulness and
validity of the proposed approach.

\medskip


\medskip

\end{abstract}


\section{Introduction}
\label{sec:intro}

The human immunodeficiency virus (HIV) is a retrovirus that infects cells of
the immune system, destroying or impairing their function. HIV is transmitted
primarily via unprotected sexual intercourse, contaminated blood transfusions,
hypodermic needles, and from mother to child during pregnancy, delivery,
or breastfeeding \cite{RomMark07}. As the infection progresses, the immune
system becomes weaker, and the person becomes more susceptible to infections.
The most advanced stage of HIV infection is acquired
immunodeficiency syndrome (AIDS)~\cite{WhoSiteHivAids}.
There is no cure or vaccine to AIDS. However, antiretroviral (ART) treatment
improves health, prolongs life, and substantially reduces the risk of HIV transmission.
In both high-income and low-income countries, the life expectancy of patients
infected with HIV who have access to ART is now measured in decades, and might
approach that of uninfected populations in patients who receive an optimum
treatment (see \cite{DeeksEtAl13} and references cited therein).
However, ART treatment still presents substantial limitations: does not fully
restore health; treatment is associated with side effects; the medications
are expensive; and is not curative.
Following the Joint United Nations Programme on HIV and AIDS (UNAIDS), in 2013
there were approximately 35 million people living with HIV globally.
An estimated 2.1 million people became newly infected with HIV in 2013,
down from 3.4 million in 2001 worldwide. The number of new HIV infection among
children has declined 58\% since 2001, being in 2013 approximately
240 000 worldwide. The number of AIDS-related deaths have fallen by 35\% since
the peak in 2005. In 2013, approximately 1.5 million people died
from AIDS-related causes worldwide. In 2013, around 12.9 million people living
with HIV had access to ART therapy, which represents, approximately, 37\%
of all people living with HIV \cite{UNAIDSFactSheet2014,UNAIDSGapRep2014}.

\emph{Mycobacterium tuberculosis} is the cause of most occurrences
of tuberculosis (TB) and is usually acquired via airborne infection
from someone who has active TB. It typically affects the lungs (pulmonary TB)
but can affect other sites as well (extrapulmonary TB).
According with the World Health Organization (WHO), in 2013,
an estimated 9.0 million people developed TB
and 1.5 million died from the disease,
360 000 of whom were HIV-positive. TB is slowly
declining each year and it is estimated that 37 million
lives were saved between 2000 and 2013 through effective
diagnosis and treatment. However, since most deaths
from TB are preventable, the death toll from the disease
is still unacceptably high and efforts to combat it must be
accelerated \cite{WHO14}.

Following WHO, the human immunodeficiency virus (HIV)
and \emph{mycobacterium tuberculosis} are the first
and second cause of death from a single infectious agent, respectively \cite{WHO13}.
Both HIV/AIDS and TB are present in all regions of the world \cite{Morison01,WHO14}.
Individuals infected with HIV are more likely to develop TB disease because of their
immunodeficiency, and HIV infection is the most powerful risk factor for progression
from TB infection to disease \cite{GetahunEtAl10}.
In 2013, 1.1 million of 9.0 million people who developed TB worldwide were HIV-positive.
The number of people dying from HIV-associated to TB has been falling since 2003. However,
there were still 360 000 deaths from HIV-associated to TB in 2013, and further efforts
are needed to reduce this burden \cite{WHO14}. ART is a critical intervention
for reducing the risk of TB morbidity and mortality among people living with HIV and,
when combined with TB preventive therapy, it can have a significant impact
on TB prevention \cite{WHO14}.
Collaborative TB/HIV activities (including HIV testing,
ART therapy and TB preventive measures) are crucial
for the reduction of TB-HIV coinfected individuals.
WHO estimates that these collaborative activities prevented
1.3 million people from dying, from 2005 to 2012. However,
significant challenges remain: the reduction of tuberculosis
related deaths among people living with HIV has slowed in recent years;
the ART therapy is not being delivered to TB-HIV coinfected patients
in the majority of the countries with the largest number of TB/HIV patients;
the pace of treatment scale-up for TB/HIV patients has slowed; less than half
of notified TB patients were tested for HIV in 2012; and only a small fraction
of TB/HIV infected individuals received TB preventive therapy \cite{UNAIDSRep2013}.
The study of the joint dynamics of TB and HIV present formidable mathematical challenges
due to the fact that the models of transmission are quite distinct \cite{ChavezEtAll09}.
Here we focus on a recent mathematical model of optimal control
for TB-HIV/AIDS coinfection proposed by \cite{SiTo15}.

Optimal control is a branch of mathematics developed to find optimal ways to control
a dynamic system \cite{PoBoGrMi62}, e.g. a dynamic system that models infectious
diseases. Optimal control has been applied to TB models, HIV models and also
co-infection models (see, e.g.,
\cite{rv2:agustu,JuLeFe02,KirsLenSer96,LeWo07,MaMuChuMu09,RoSiTo14,SiTo13,SiTo15}
and references cited therein for TB-HIV/AIDS models and \cite{rv2:Okosun} 
for co-infection of malaria and cholera). 
In this paper we consider the optimal control
problem for the TB-HIV/AIDS model proposed in \cite{SiTo15} from a multiobjective
perspective. Our approach avoids the use of weight parameters and allows
to obtain a wide range of optimal control strategies. These strategies offer
the decision maker useful information for effective decision making.

Traditional mathematical programming methods for solving multiobjective
optimization problems (MOPs) convert the original problem into a single-objective
optimization problem. This is referred as to scalarization and the function
to be optimized, which depends on some parameters, is termed the scalarizing function.
A solution to the scalarizing function, obtained using a single-objective
optimization algorithm, is expected to be Pareto optimal. For approximating
multiple Pareto optimal solutions, repeated runs with different
parameter settings must be performed.
The weighted sum method~\cite{GaSa55} consists in minimizing a weighted sum of
multiple objectives. For problems with a convex Pareto front, this method
guarantees finding solutions in the entire Pareto optimal region. However,
it fails to find solutions in nonconvex regions of the Pareto front. Weighted
metric methods~\cite{Mi99} are based on minimizing a weighted distance between
some reference point and the feasible objective region. The widely used approach
belonging to this class of methods is the Chebyshev method~\cite{Bo76},
which consists in minimizing a weighted infinity norm. Although solutions
in convex and nonconvex regions of the Pareto front can be obtained by this method,
a resulting scalarizing function becomes nondifferentiable even when all the
objectives are differentiable. The problem resulting from the Chebyshev method
can be reformulated in the smooth form. The resulting formulation is known
as the goal attainment method~\cite{Mi99} or the Pascoletti--Serafini
scalarization~\cite{PaSe84}. In this method, a slack variable is minimized
and the weighted difference for each objective is converted into a constraint.
Although the problem can be solved in a differentiable form, problem complexity
is augmented by adding one additional variable and $m$ constraints (where $m$
is the number of objectives). The normal boundary intersection and normal constraint
methods use a hyperplane with uniformly distributed points passing through
the critical points of the Pareto front. The normal boundary intersection
method~\cite{DaDe98} searches for the maximum distance from a point on the
simplex along the normal pointing toward the origin. The obtained point may
or may not be a Pareto optimal point, with the resulting scalarizing problem
including an equality constraint that is not easy to treat for all the cases.
On the other hand, the normal constraint method~\cite{Me04} uses an inequality
constraint reduction of the feasible objective space and the normalized function
values to cope with disparate function scales. The method is successful in achieving
a uniform distribution of approximating points, though there is no guarantee
that an obtained point is Pareto optimal. Here, motivated by the results
recently obtained in~\cite{DeSiTo15} for a TB model and in \cite{MyID:325,MyID:333}
for the dengue disease, we adopt the $\epsilon$-constraint method~\cite{HaLaWi71}.
This method suggests optimizing one objective function
and converting all other objectives as constraints, setting an upper bound
to each of them. Solutions obtained using multiobjective optimization
provide comprehensive insights about the optimal strategies and the diseases
dynamics resulting from implementation of those strategies.

The paper is organized as follows. In Section~\ref{sec:model} we briefly
describe the TB-HIV/AIDS model from \cite{SiTo15}. The multiobjective optimization theory
is applied to an optimal control problem in Section~\ref{sec:optContolProb}:
we start by formulating the optimal control problem in Subsection~\ref{subsecOCprob},
then we consider this problem from a multiobjective perspective
(Subsection~\ref{subsecMultiobj}) and we describe the numerical method that we
use to solve the multiobjective problem (Subsection~\ref{subsec:scalar}).
In Section~\ref{sec:scalar} we present and discuss numerical results
for the multiobjective problem. We end with Section~\ref{sec:conc}
of conclusions and future work.


\section{TB-HIV/AIDS coinfection model}
\label{sec:model}

The present study considers the population model for TB-HIV/AIDS coinfection
proposed in~\cite{SiTo15}, where TB, HIV and TB-HIV infected individuals
have access to respective disease treatment, and single HIV-infected
and TB-HIV co-infected individuals under HIV and TB/HIV treatment,
respectively, stay in a \emph{chronic} stage of the HIV infection. The model
divides the population into eleven mutually exclusive compartments: susceptible
individuals ($S$); TB-latently infected individuals, who have no symptoms
of TB disease and are not infectious ($L_T$); TB-infected individuals,
who have active TB disease and are infectious ($I_T$); TB-recovered individuals
($R$); HIV-infected individuals with no clinical symptoms of AIDS ($I_H$);
HIV-infected individuals under treatment for HIV infection ($C_H$);
HIV-infected individuals with AIDS clinical symptoms ($A$); TB-latent individuals
co-infected with HIV (pre-AIDS) ($L_{TH}$); HIV-infected individuals (pre-AIDS)
co-infected with active TB disease ($I_{TH}$); TB-recovered individuals with
HIV-infection without AIDS symptoms ($R_{H}$); and HIV-infected individuals
with AIDS symptoms co-infected with active TB ($A_T$).
The total population at time $t$, denoted by $N(t)$, is given by
\begin{multline*}
N(t) = S(t) + L_T(t) + I_T(t) + R(t) + I_H(t) + A(t) + C_H(t) \\
+ L_{TH}(t) + I_{TH}(t) + R_H(t) + A_T(t).
\end{multline*}
Susceptible individuals acquire TB infection from individuals
with active TB at a rate $\lambda_T$, given by
\[
\lambda_T(t) = \frac{\beta_1}{N(t)} (I_T(t) + I_{TH}(t) + A_T(t)),
\]
where $\beta_1$ is the effective contact rate for TB infection. Similarly,
susceptible individuals acquire HIV infection, following effective contact
with people infected with HIV at a rate $\lambda_H$, given by
\[
\lambda_H(t) = \frac{\beta_2}{N(t)} [I_H(t) + I_{TH}(t) + L_{TH}(t)
+ R_H(t) + \eta_C C_H(t) + \eta_A (A(t) + A_T(t))],
\]
where $\beta_2$ is the effective contact rate for HIV transmission,
$\eta_A \geq 1$ is the modification parameter that accounts for the relative
infectiousness of individuals with AIDS symptoms and $\eta_C \leq 1$ is the
modification parameter that translates the partial restoration of immune
function of individuals with HIV infection that use correctly the antiretroviral
treatment. The remaining parameters used to describe
the model are presented in Table~\ref{tab:params}.
\begin{table}[t]
\centering
\small
\begin{tabular}{lll}
\toprule
Symbol & Description  & Value \\
\hline
$T$ & Considered time in years & 10 \\
$N(0)$ & Initial population size & 30000 \\
$\gamma_1$ & Modification parameter & 0.9 \\
$\gamma_2$ & Modification parameter & 1.1 \\
$\eta_C$ & Modification parameter & 0.9 \\
$\eta_A$ & Modification parameter & 1.05 \\
$\delta$ & Modification parameter & 1.03 \\
$\psi$ & Modification parameter & 1.07 \\
$\beta_1$ & TB transmission rate & 0.6 \\
$\beta_2$ & HIV transmission rate & 0.1 \\
$\mu$ & Recruitment rate & 430.0 \\
$k_1$ & Rate at which individuals leave $L_T$ class by becoming infectious & 1.0/2.0 \\
$k_2$ & Rate at which individuals leave $L_{TH}$ class by becoming TB infectious & $1.3 k_1$ \\
$k_3$ & Rate at which individuals leave $L_{TH}$ class & 2.0 \\
$\rho_1$ & Rate at which individuals leave $I_H$ class to $A$ & 0.1 \\
$\rho_2$ & Rate at which individuals leave $I_{TH}$ class & 1.0 \\
$\omega_1$ & Rate at which individuals leave $C_H$ class & 0.09 \\
$\omega_2$ & Rate at which individuals leave $R_H$ class & 0.15 \\
$\tau_1$ & TB treatment rate for $L_T$ individuals & 2.0 \\
$\tau_2$ & TB treatment rate for $I_T$ individuals & 1.0 \\
$\phi$ & HIV treatment rate for $I_H$ individuals & 1.0 \\
$\alpha_1$ & AIDS treatment rate & 0.33 \\
$\alpha_2$ & HIV treatment rate for $A_T$ individuals & 0.33 \\
$r$ & Fraction of $L_{TH}$ individuals that take HIV and TB treatment & 0.3 \\
$d_N$ & Natural death rate & 1.0/70.0 \\
$d_T$ & TB induced death rate & 1.0/10.0 \\
$d_A$ & AIDS induced death rate & 0.3 \\
$d_{TA}$ & AIDS-TB induced death rate & 0.33 \\
\bottomrule
\end{tabular}
\caption{Model parameters, borrowed from~\cite{SiTo15}.}\label{tab:params}
\end{table}

Two control functions, which represent prevention and treatment measures, are
introduced into the model and are continuously implemented during a considered
period of disease treatment: the control $u_1(t)$ represents the fraction
of $I_{TH}$ individuals that takes HIV and TB treatment, simultaneously;
$u_2(t)$ represents the fraction of $I_{TH}$ individuals
that takes TB treatment only  \cite{SiTo15}.

\begin{figure}[!htb]
	\centering
	\includegraphics[scale=0.7]{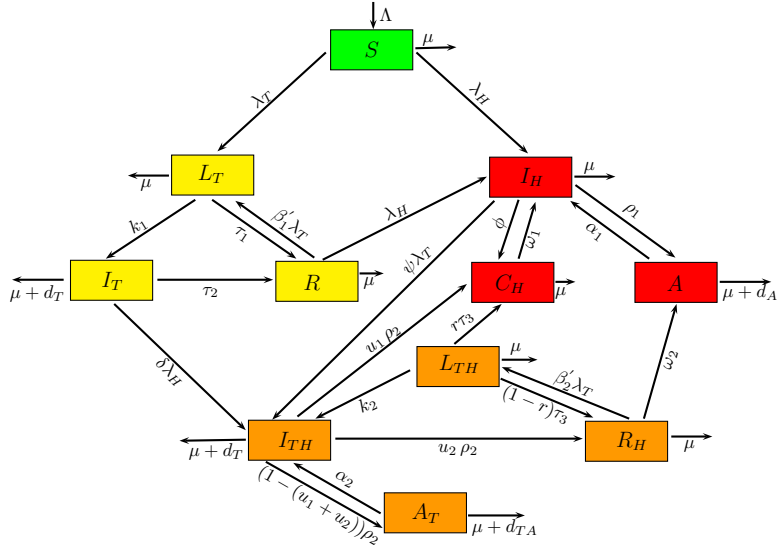}
	\caption{Model for TB-HIV/AIDS transmission.}
	\label{fig:model:flow}
\end{figure}

The transmission dynamics of TB-HIV/AIDS coinfection is modeled
by the following system of differential equations:
\begin{equation}
\label{eq:model}
\small
\left\{
\begin{array}{l}
\dot{S}(t) = \mu - \lambda_T(t) S(t) - \lambda_H(t) S(t) - d_N S(t),\\
\dot{L}_T(t) = \lambda_T(t) S(t) + \gamma_1 \lambda_T(t) R(t) - (k_1 + \tau_1 + d_N) L_T(t),\\
\dot{I}_T(t) = k_1 L_T(t) - (\tau_2 + d_T + d_N + \delta \lambda_H(t)) I_T(t),\\
\dot{R}(t) = \tau_1 L_T(t) + \tau_2 I_T(t) - (\gamma_1 \lambda_T(t)
+ \lambda_H(t) + d_N) R(t),\\
\dot{I}_H(t) = \lambda_H(t) S(t) - (\rho_1 + \phi + \psi \lambda_T(t) + d_N) I_H(t)
+ \alpha_1 A(t) + \lambda_H(t) R(t) + \omega_1 C_H(t),\\
\dot{A}(t) = \rho_1 I_H(t) + \omega_2 R_H(t) - \alpha_1 A(t) - (d_N + d_A) A(t),\\
\dot{C}_H(t) = \phi I_H(t) + u_1(t) \rho_2 I_{TH}(t) + r k_3 L_{TH}(t)
- (\omega_1 + d_N) C_H(t),\\
\dot{L}_{TH}(t) = \gamma_2 \lambda_T(t) R_H(t) - (k_2 + k_3 + d_N) L_{TH}(t),\\
\dot{I}_{TH}(t) = \delta \lambda_H(t) I_T(t) + \psi \lambda_T(t) I_H(t)
+ \alpha_2 A_T(t) + k_2 L_{TH}(t) - (\rho_2 + d_N + d_T) I_{TH}(t),\\
\dot{R}_H(t) = u_2(t) \rho_2 I_{TH}(t) + (1 - r) k_3 L_{TH}(t)
- (\gamma_2 \lambda_T(t) + \omega_2 + d_N) R_H(t),\\
\dot{A}_T(t) = (1 - (u_1(t) + u_2(t))) \rho_2 I_{TH}(t)
- (\alpha_2 + d_N + d_{TA}) A_T(t).
\end{array}
\right.
\end{equation}
The model flow is illustrated in Figure~\ref{fig:model:flow}.
The initial conditions are given in Table~\ref{tab:init:conds}.
\begin{table}[t]
\centering
\small
\begin{tabular}{lll} 
\toprule 
Cathegory & Description  & Initial Value \\[0.1cm]
\hline
$S$ & Susceptible & $\frac{66N(0)}{120}$ \\[0.1cm]
${L}_T$ & TB-Latent & $\frac{37N(0)}{120}$\\[0.1cm]
$I_T$ & TB-Active infected  & $\frac{5N(0)}{120}$\\[0.1cm]
$R$ & TB-Recovered & $\frac{37N(0)}{120}$ \\[0.1cm]
$I_H$ & HIV-Infected (pre-AIDS)& $\frac{2N(0)}{120}$ \\[0.1cm]
$A$ & HIV-infected with AIDS symptoms & $\frac{37N(0)}{120}$ \\[0.1cm]
$C_H$ & HIV-infected under ART therapy & $\frac{N(0)}{120}$ \\[0.1cm]
$L_{TH}$ & TB-Latent co-infected with HIV (pre-AIDS) & $\frac{2N(0)}{120}$ \\[0.1cm]
$I_{TH}$ & HIV-Infected (pre-AIDS) co-infected with active TB & $\frac{2N(0)}{120}$ \\[0.1cm]
$R_H$ & TB-recovered with HIV-infection (pre-AIDS) & $\frac{N(0)}{120}$ \\[0.1cm]
$A_T$ & HIV-Infected with AIDS symptoms co-infected with active TB & $\frac{N(0)}{120}$ \\[0.1cm]
\bottomrule 
\end{tabular}
\caption{Initial conditions for the state variables of the TB-HIV/AIDS
model \eqref{eq:model}, borrowed from~\cite{SiTo15}.}\label{tab:init:conds}
\end{table}
Note that the period spent in class $I_{TH}$ does not change with the control measures
because the controls represent the fraction of individuals that are treated
both for TB and HIV and only for TB, and not the treatment duration. Indeed, 
the period spent in class $I_{TH}$ is given by the constant treatment rate $\rho_2$.


\section{Multiobjective approach to an optimal control problem}
\label{sec:optContolProb}

Traditionally, the problem of finding a control law for a given system
is addressed by optimal control theory \cite{PoBoGrMi62}.


\subsection{Optimal control problem}
\label{subsecOCprob}

In the optimal control approach, the aim is to find the optimal values
$u_1^*$ and $u_2^*$ of the controls $u_1$ and $u_2$, such that the associated
state trajectories $S^*$, $L_T^*$, $I_T^*$, $R^*$, $I_H^*$, $A^*$, $C_H^*$,
$L_{TH}^*$, $I_{TH}^*$, $R_{H}^*$, $A_{T}^*$, are solution of
system~\eqref{eq:model} in the time interval $[0, T]$, with the initial
conditions in Table~\ref{tab:init:conds}, and minimize an objective functional.
Consider the state system of ordinary differential equations~\eqref{eq:model}
and the set of admissible control functions given by
\begin{align*}
\Omega = \{ (u_1(\cdot),u_2(\cdot))
\in (L^{\infty}(0,T))^2 \, | \, & 0 \leq u_1(t),u_2(t) \leq 0.95 \\
& \, \wedge \, u_1(t) + u_2(t) \leq 0.95, \, \forall t \in [0,T] \}.
\end{align*}
According to~\cite{SiTo15}, the objective functional can be defined as
\begin{equation}
\label{oc:problem}
J(u_1(\cdot),u_2(\cdot))
= \int_{0}^{T}\left[ A(t) + A_T(t) + w_1u_1^2(t) + w_2u_2^2(t) \right]dt,
\end{equation}
where the constants $w_1$ and $w_2$ are a measure of the relative cost of the
interventions associated with the controls $u_1$ and $u_2$, respectively.
Note that the objective functional \eqref{oc:problem} is a function of state
and control variables. Its minimization implies three important aspects:
(i) reducing the number of individuals with AIDS symptoms, (ii) decreasing
the number of individuals with AIDS symptoms and active TB disease and
(iii) reducing the costs of implementing treatment policies. The optimal
control problem consists in determining ($S^*$, $L_T^*$, $I_T^*$, $R^*$,
$I_H^*$, $A^*$, $C_H^*$, $L_{TH}^*$, $I_{TH}^*$, $R_{H}^*$, $A_{T}^*$),
associated to admissible controls
$\left(u_1^*(\cdot), u_2^*(\cdot)\right) \in \Omega$
on the time interval $[0, T]$, satisfying~\eqref{eq:model}, the initial
conditions in Table~\ref{tab:init:conds}, and minimizing the objective
functional~\eqref{oc:problem}, \textrm{i.e.},
\[
J(u_1^*(\cdot), u_2^*(\cdot)) = \min_{\Omega} J(u_1(\cdot), u_2(\cdot)).
\]
The approach based on optimal control theory adopted in \cite{SiTo15}
allows to obtain the optimal solution to the cost functional~\eqref{oc:problem},
which is defined from some decision maker's perspective by means of the constants
$w_1$ and $w_2$. However, the choice of the values of $w_1$ and $w_2$
requires some knowledge about the problem and the decision maker's preferences,
which often are not available in advance. Another disadvantage consists in the
fact that a single optimal solution to~\eqref{oc:problem}
does not provide all useful insights about the optimal strategies
and corresponding dynamics. A large range of alternatives remain unexplored
and the decision maker is limited in his/her options. In~\cite{SiTo15}, numerical
simulations to the optimal control problem are performed using $(w_1, w_2)=(25,25)$
and $(w_1, w_2)=(250,25)$. Both values are larger than one, which suggests that
they are adapted to the scale of the objectives. The choice of the values
of these parameters is not straightforward and there is no guarantee
that the best compromise solution has been obtained.


\subsection{Multiobjective optimization}
\label{subsecMultiobj}

Our work addresses the optimal control problem for the TB-HIV/AIDS
coinfection model~\eqref{eq:model} from a multiobjective perspective.
A multiojective optimization problem is formulated by decomposing
the cost functional shown in~\eqref{oc:problem} into two components,
representing different aspects that must be taken into consideration
when dealing with TB-HIV/AIDS. The problem of finding
the optimal controls is defined as:
\begin{equation}
\label{mo:problem}
\begin{array}{rl}
\text{minimize}
& f_1(A(\cdot),A_T(\cdot))
= \displaystyle \int_{0}^{T} \left[ A(t) + A_T(t) \right] dt, \\[0.2cm]
& f_2(u_1(\cdot),u_2(\cdot))
= \displaystyle \int_{0}^{T} \left[ u_1^2(t) + u_2^2(t) \right] dt, \\[0.2cm]
\text{subject to}
& 0 \leq u_1(t) \leq 0.95, \\
& 0 \leq u_2(t) \leq 0.95, \\
& u_1(t) + u_2(t) \leq 0.95.
\end{array}
\end{equation}
In the above formulation, the weights are absent and the two objectives represent
the medical and economical perspectives, respectively. This naturally reflects
the conflicting nature of the underlying decision making problem, hence,
solving \eqref{mo:problem} is interesting and challenging.


\subsection{Scalarization}
\label{subsec:scalar}

A traditional mathematical programming approach to solving
a multiojective optimization problem consists in transforming 
an original problem with multiple objectives into a number 
of single-objective subproblems. This is referred to 
as \emph{scalarization}. The transformation is performed 
by means of a scalarizing function with some user-defined 
parameters. A single Pareto optimal solution is sought 
by optimizing each subproblem. Repeated runs with different 
parameter settings for the scalarizing function are used 
to approximate multiple Pareto optimal solutions. 

Several approaches to scalarization have been developed. 
They differ in the way the scalarizing function is formulated. 
The weighted sum method suggests minimizing a weighted sum 
of the objectives~\cite{GaSa55}. The limitation of this method 
is that solutions can only be obtained in convex regions 
of the Pareto front. On the other hand, the $\epsilon$-constraint 
method~\cite{HaLaWi71} suggests optimizing one objective function
and converting all other objectives into constraints by setting 
an upper bound to each of them. This method can find solutions 
in both convex and nonconvex regions of the Pareto front.
The method of weighted metrics~\cite{Mi99} seeks to minimize 
the distance between the feasible objective region and some 
reference point. This method is also known as 
\emph{compromise programming}~\cite{Zeleny1976}. For measuring 
the distance, a weighted $L_p$ norm is utilized. 
When the value of $p$ is small, the method may fail 
to find solutions in nonconvex regions. When $p=\infty$, 
the method defines the weighted Chebyshev problem~\cite{Bo76}.
This problem consists in minimizing the largest weighted 
deviation of one objective. By optimizing the weighted 
Chebyshev problem, solutions from convex and nonconvex 
regions can be generated. A major drawback is that 
even when the original MOP is differentiable, 
the resulting single-objective problem is nondifferentiable.
Weakly Pareto optimal solutions can be also obtained~\cite{Mi99}. 
A relaxed formulation of the Chebyshev problem with differentiable 
scalarizing function is known as the Pascoletti--Serafini 
scalarization~\cite{PaSe84}. Though, this method introduces 
one additional variable and one constraint for each objective function. 
The normal boundary intersection~\cite{DaDe98} uses a hyperplane 
with evenly distributed points that passes through the extreme points 
of the Pareto front. For each point on the hyperplane, the method 
searches for the maximum distance along the normal pointing 
toward the origin. The normal constraint method~\cite{Me04} 
suggests optimizing one objective and employing an inequality 
constraint reduction of the feasible space using the points 
on the hyperplane. However, there is no guarantee that the solutions 
obtained by the normal boundary intersection and normal constraint methods 
are Pareto optimal. A comparative analysis of the different scalarization 
approaches on optimal control problems from epidemiology 
can be found in~\cite{MyID:333,DeSiTo15}.

Motivated by the results recently obtained in~\cite{DeSiTo15} for a TB model,
we adopt here the $\epsilon$-constraint method~\cite{HaLaWi71}. This method
suggests optimizing one objective function while converting all other objectives
into constraints by setting an upper bound to each of them. The problem to be solved
is of the following form:
\begin{equation}
\label{method:eps}
\begin{array}{rlll}
\underset{\boldsymbol{x} \in \Omega}{\text{minimize}} & f_l(\boldsymbol{x}) & \\
\text{subject to} & f_i(\boldsymbol{x})\leq \epsilon_i,
& \forall i \in \{1,\ldots,m\} \wedge i \neq l,
\end{array}
\end{equation}
where the $l$th objective is minimized, the parameter $\epsilon_i$ represents
an upper bound of the value of $f_i$ and $m$ is the number of objectives. The
major reasons for using this method are as follows. The $\epsilon$-constraint
method is able to find solutions in convex and nonconvex regions of the Pareto
optimal front. When all the objective functions in the MOP are convex, problem
\eqref{method:eps} is also convex and has a unique solution. For any given
upper bound vector $\boldsymbol{\epsilon} = \{\epsilon_1, \ldots, \epsilon_{m-1}\}$,
the unique solution of problem~\eqref{method:eps} is Pareto optimal~\cite{Mi99}.
Moreover, when considering different scenarios in the model, the optimal solutions
obtained for the same values of $\epsilon$ can be used for comparison, as they
will lie on the same line in the objective space determined by the corresponding
value of $\epsilon$. This characteristic is convenient and helpful for the
analysis of the dynamics in the TB-HIV/AIDS model.


\section{Numerical simulations and discussion}
\label{sec:scalar}

This section presents and discusses numerical results for the optimal
controls using the multiobjective optimization approach. Moreover,
possible scenarios of applying the control strategies are investigated.


\subsection{Experimental setup}
\label{sec:setup}

The fourth-order Runge--Kutta method is used for numerically integrating
system \eqref{eq:model}. The control and state variables are discretized
using $100$ equally spaced time intervals over the period $[0,T]$.
The integrals defining objective functions in~\eqref{mo:problem}
are calculated using the trapezoidal rule.

Using the formulation~\eqref{method:eps}, the first objective
in~\eqref{mo:problem} is minimized and the second objective
is set as the constraint bounded by the value of $\epsilon$.
The Pareto front is discretized by defining 100 evenly distributed
values  of $\epsilon$ over the range of the second objective, $f_2$,
which can be calculated as $f_2^{\min}$ for $\left(u_1(\cdot), u_2(\cdot)\right)
\equiv 0$ and $f_2^{\max}$ for $u_1(t) + u_2(t) \equiv 0.95$.
It is worth noting that the $\epsilon$-constraint method does not need
the information about the range of $f_1$, which is not known beforehand
due to the presence of the constraint imposed on $u_1(t) + u_2(t)$,
whereas the majority of methods discussed in Section~\ref{sec:scalar}
may need this information.
To solve the problems with different values of $\epsilon$,
the MATLAB\textsuperscript{\textregistered} function \texttt{fmincon}
with a sequential quadratic programming algorithm is used, setting
the maximum number of function evaluations to $50,000$.
The behavior of the dynamics in the TB-HIV/AIDS model are investigated
for the cases when both the controls $u_1$ and $u_2$ are applied separately
and simultaneously. In the following, the notation for resulting multiobjective
optimization problems (MOPs) will be as shown in Table~\ref{tab:mops}:
MOP1 refers to the case when $u_1$ and $u_2$ are applied simultaneously;
MOP2 refers to the case when $u_1$ is applied alone; MOP3 refers to the case
when $u_2$ is applied alone. The components of the vector
$\boldsymbol{f} = (f_1, f_2)^{\text{T}}$ are as shown in~\eqref{mo:problem}.
\begin{table}[H]
\centering
\begin{tabular}{ccc}
\toprule
MOP1 & MOP2 & MOP3 \\
\hline
$\boldsymbol{f}(u_1(\cdot), u_2(\cdot))$
& $\boldsymbol{f}(u_1(\cdot))$ & $\boldsymbol{f}(u_2(\cdot))$ \\
\bottomrule
\end{tabular}
\caption{Different multiobjective optimization problems.}\label{tab:mops}
\end{table}


\subsection{Results and discussion}

\begin{figure}
\centering
\includegraphics[width = 0.7\textwidth]{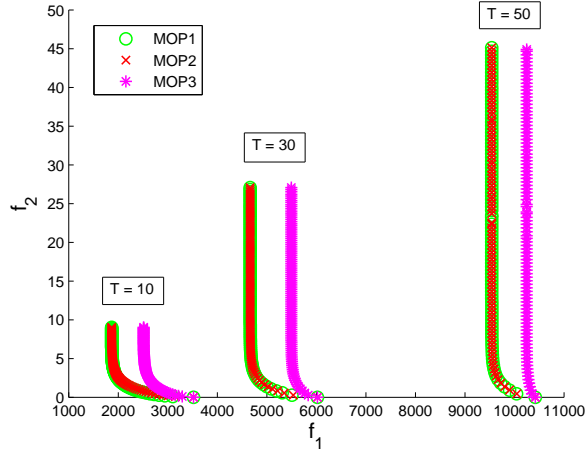}
\caption{Trade-off curves for various periods of time.}\label{fig:pf:T}
\end{figure}

Our study investigates optimal control strategies for different
treatment periods. The trade-off curves for the treatment scenarios
presented in Table~\ref{tab:mops}, considering $T \in \{10, 30, 50\}$,
are displayed in Figure~\ref{fig:pf:T}. Overall, the larger period of study,
the higher the number of individuals with AIDS and active TB. As expected,
reducing the amount of control measures leads to the increase in the number
of individuals with AIDS and active TB, whereas decrease in the number
of individuals with AIDS and active TB can be achieved though rising expenses
for treatment. These results clearly reflect the conflicting nature of the
two objectives. Also, it can be seen that an efficient range of the control
policies is limited, as starting from some point the reduction in the number
of individuals with AIDS and active TB is possible through exponential increase
in expenses for medication. Since available resources are often scarce, scenarios
involving high expenses may be practically unacceptable. All curves share
similar features and similar trends can be identified through analysis of the
obtained solutions. Due to these facts and space limitation, in what follows
the obtained results are discussed only for $T=10$.

\begin{figure}[!hb]
\centering
\includegraphics[width = 0.7\textwidth]{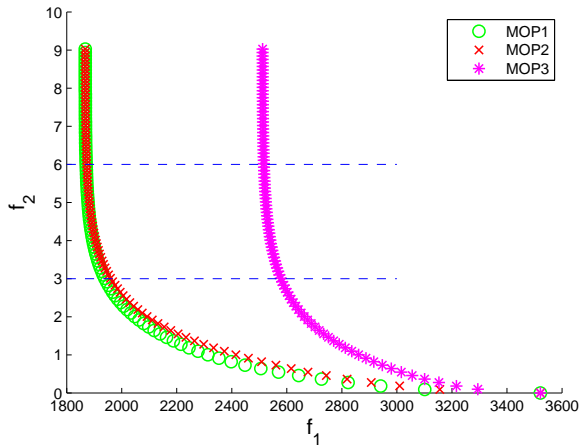}
\caption{Trade-off curves for $T=10$.}\label{fig:pf}
\end{figure}

Figure~\ref{fig:pf} shows the trade-off curves obtained for MOP1, MOP2 and MOP3,
corresponding to $T=10$. As one can see, the three curves share a common point
in the objective space. This represents an economic perspective, i.e., there is
no treatment of TB and HIV/AIDS, the only focus is on saving money. Naturally,
this leads to an uncontrollable spread of the diseases and higher numbers
of individuals with AIDS and active TB. As control policies start being implemented,
the response of $(A+A_T)$ differs for the three considered cases. The best scenarios
from the medical perspective, i.e., when the maximum amounts of controls are applied,
are different. Interestingly, the lowest number of $(A+A_T)$ is achieved when only
implementing $u_1(\cdot)\equiv 0.95$, which corresponds to the treatment of patients
for HIV/AIDS and TB together. In this case, scenarios resulting from MOP1 and MOP2
are identical, being represented by the same point in the objective space. However,
when optimal strategies involve the treatment for TB, this allows to decrease
the number of $(A+A_T)$ in scenarios representing trade-off between the economic
and medical perspectives. It can be understood observing all the intermediate
solutions for MOP1 in Figure~\ref{fig:pf}, which give a less value of $f_1$
when comparing with solutions for MOP2 involving the same amount of control
measures. Though, treating patients only for TB appears to be an ineffective
approach when comparing with scenarios represented by MOP1 and MOP2, as optimal
solutions give significantly larger numbers of $(A+A_T)$
along the whole Pareto optimal region.
\begin{figure}
\centering
\includegraphics[width = 0.4\textwidth]{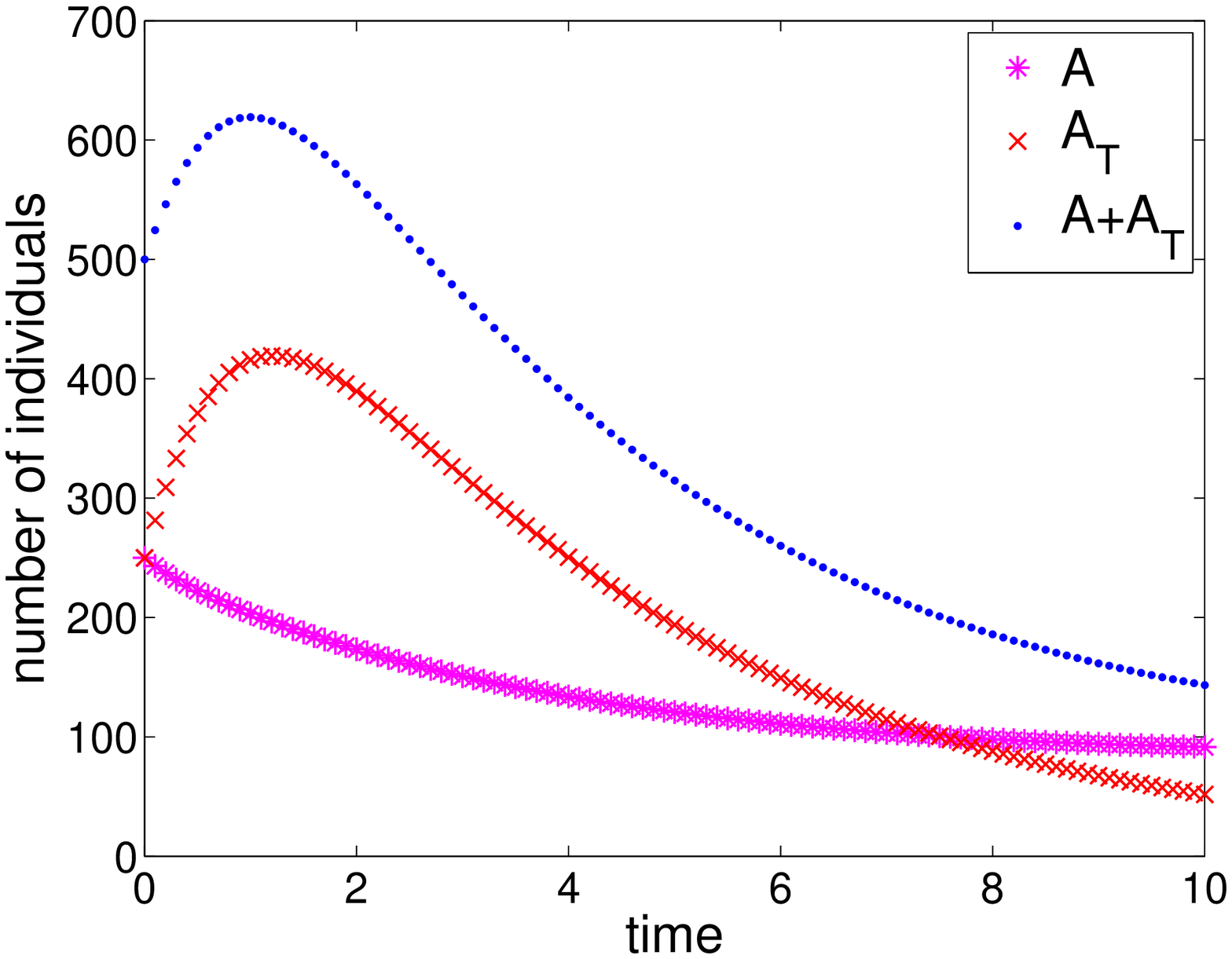}
\caption{Dynamics of individuals with AIDS and active TB without controls.}\label{fig:f0:inds}
\centering
\subfigure[MOP1 and MOP2]{\epsfig{file=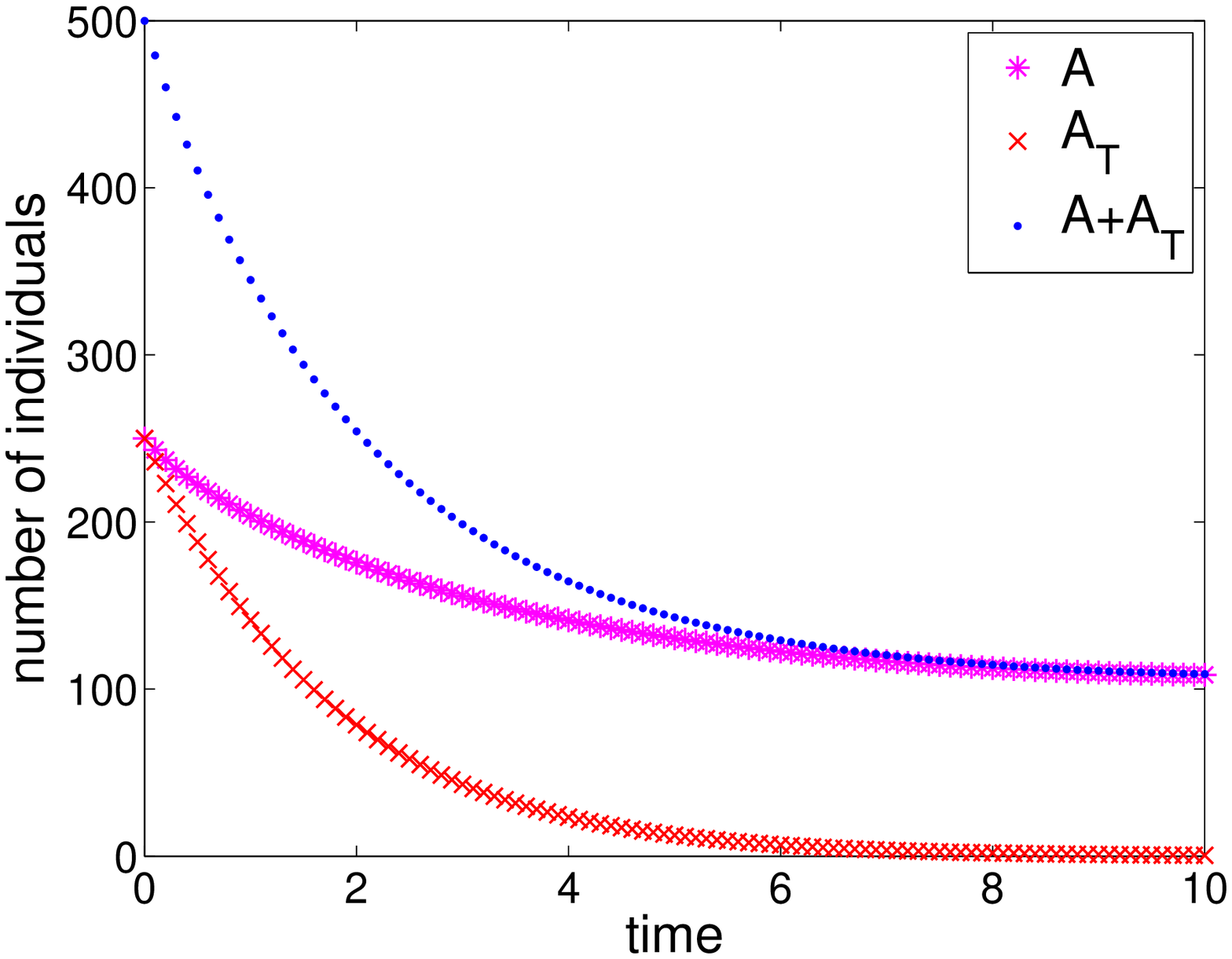,width
= 0.33\textwidth}\label{fig:max:control:MOP12}}%
\subfigure[MOP3]{\epsfig{file=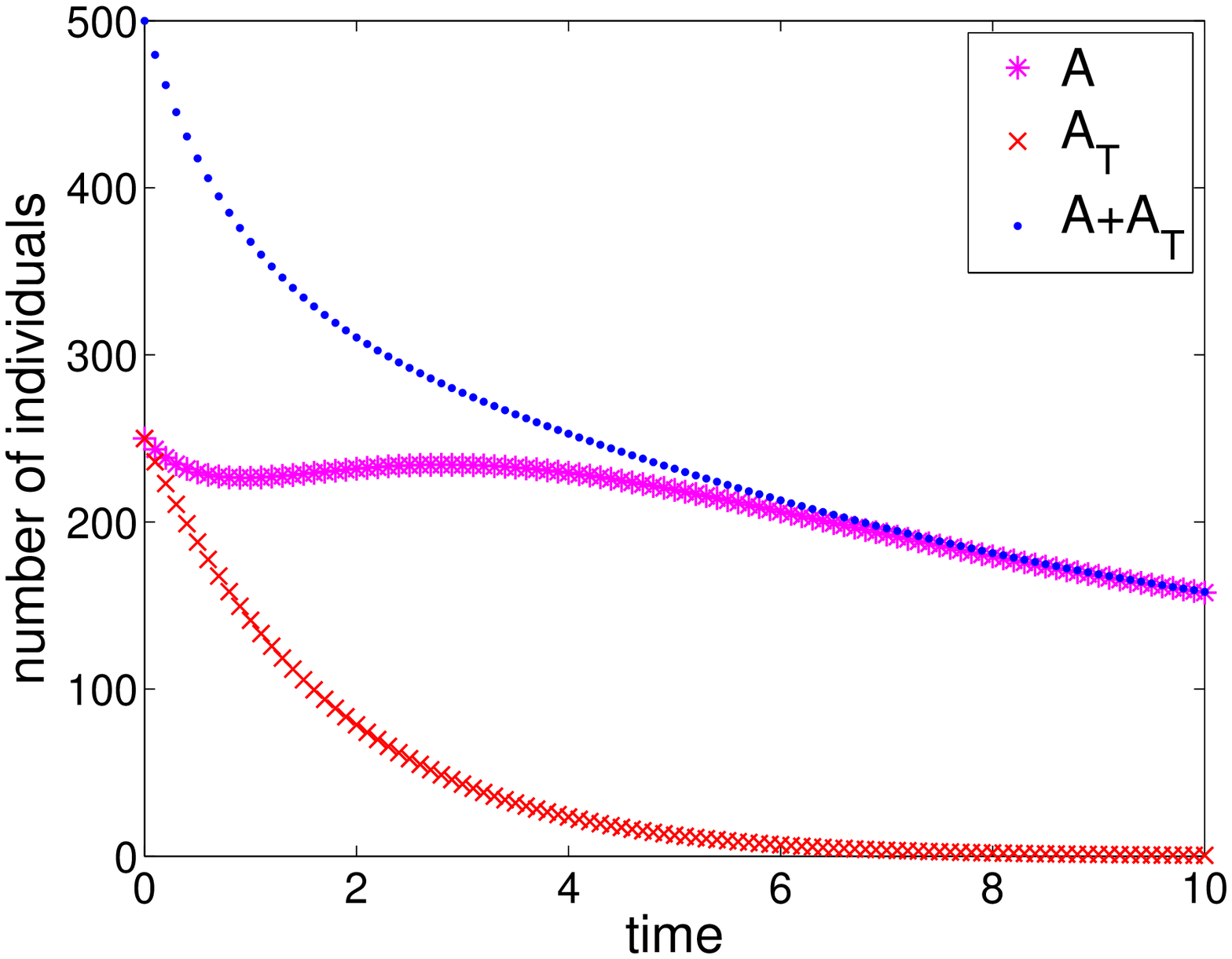,width = 0.33\textwidth}\label{fig:max:control:MOP3}}%
\caption{Dynamics of individuals with AIDS and active TB with maximum controls.}\label{fig:max:control}
\end{figure}
\begin{figure}
\centering
\subfigure[MOP1]{\epsfig{file=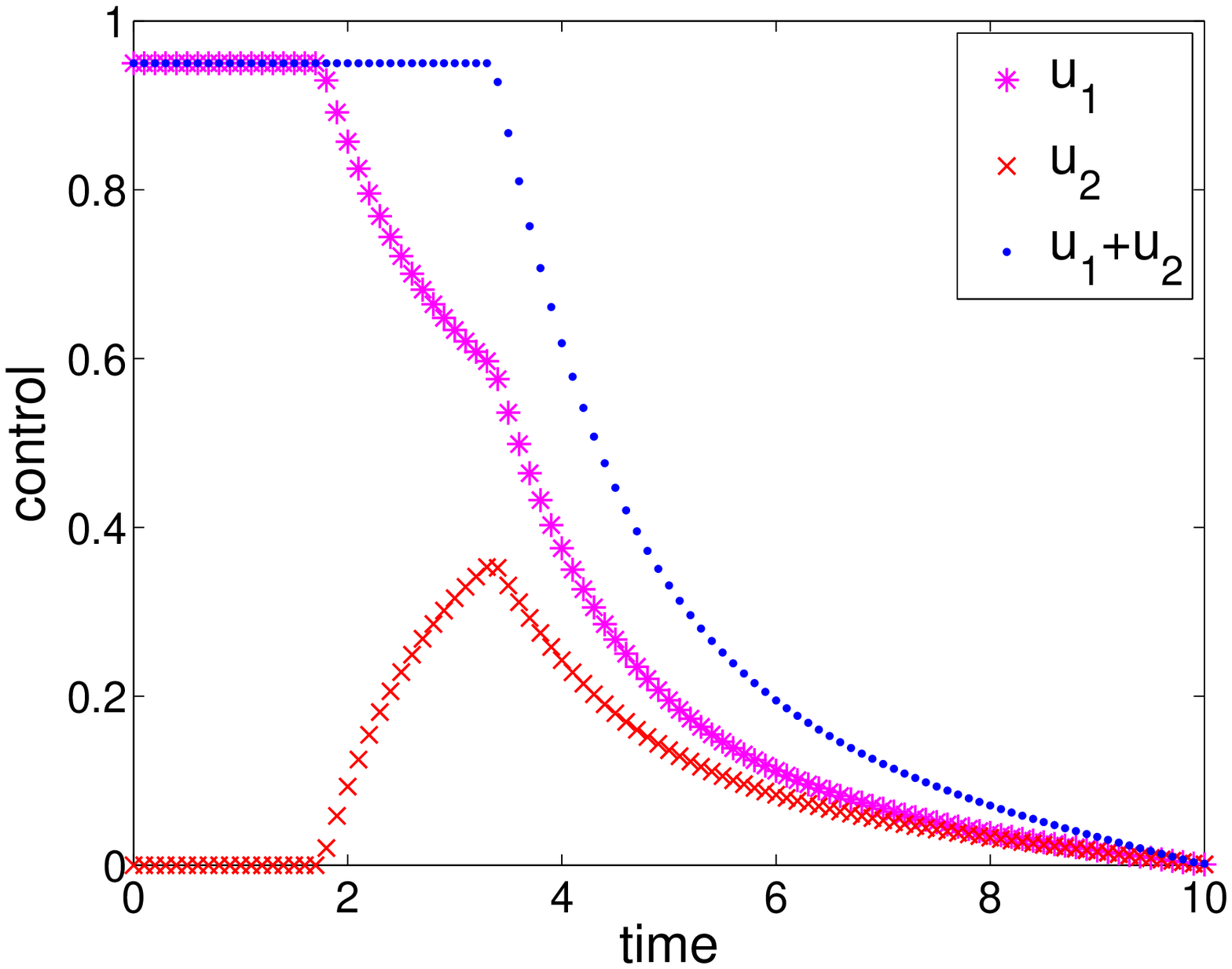,width = 0.33\textwidth}\label{fig:f3:control:MOP1}}%
\subfigure[MOP2]{\epsfig{file=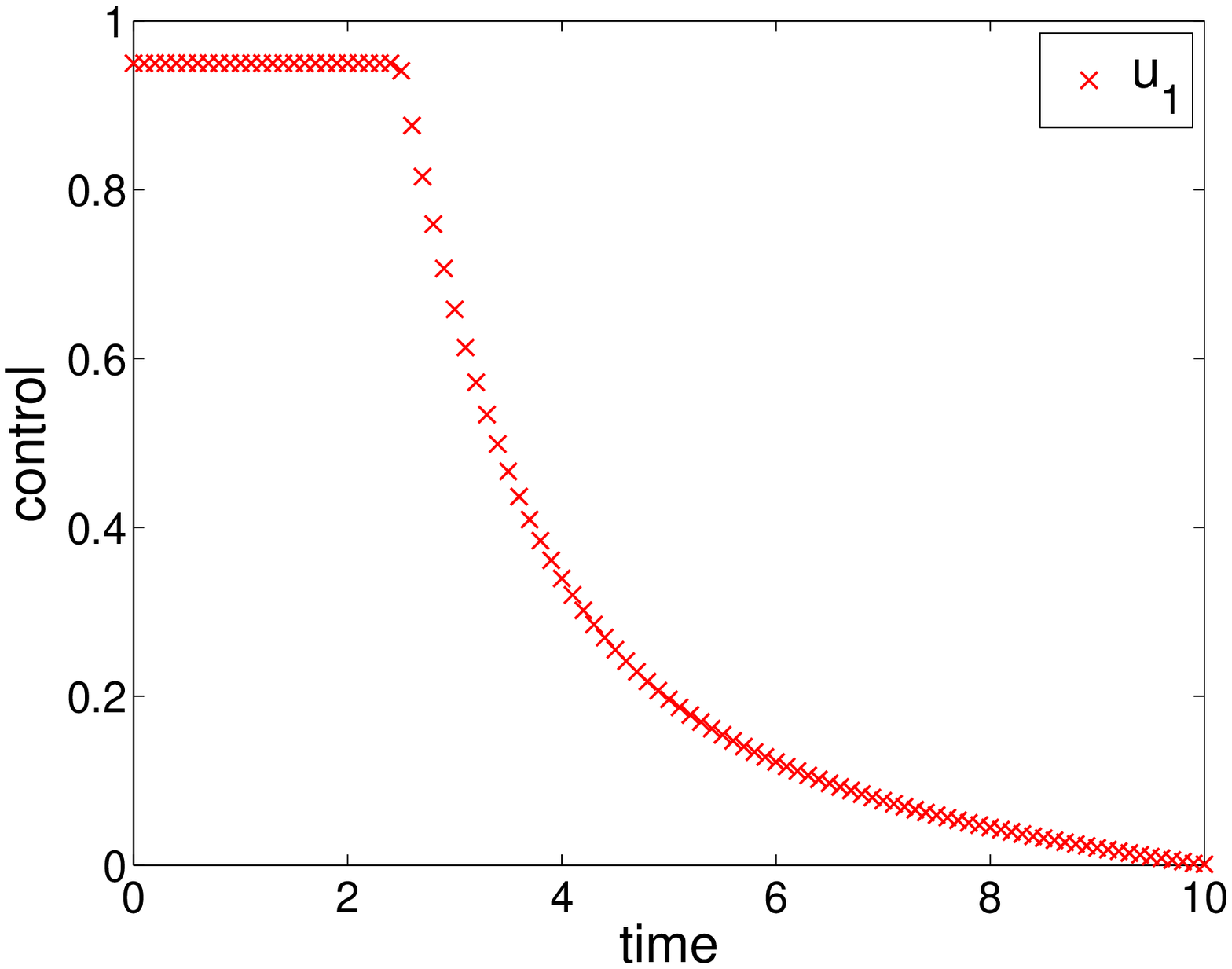,width = 0.33\textwidth}\label{fig:f3:control:MOP2}}%
\subfigure[MOP3]{\epsfig{file=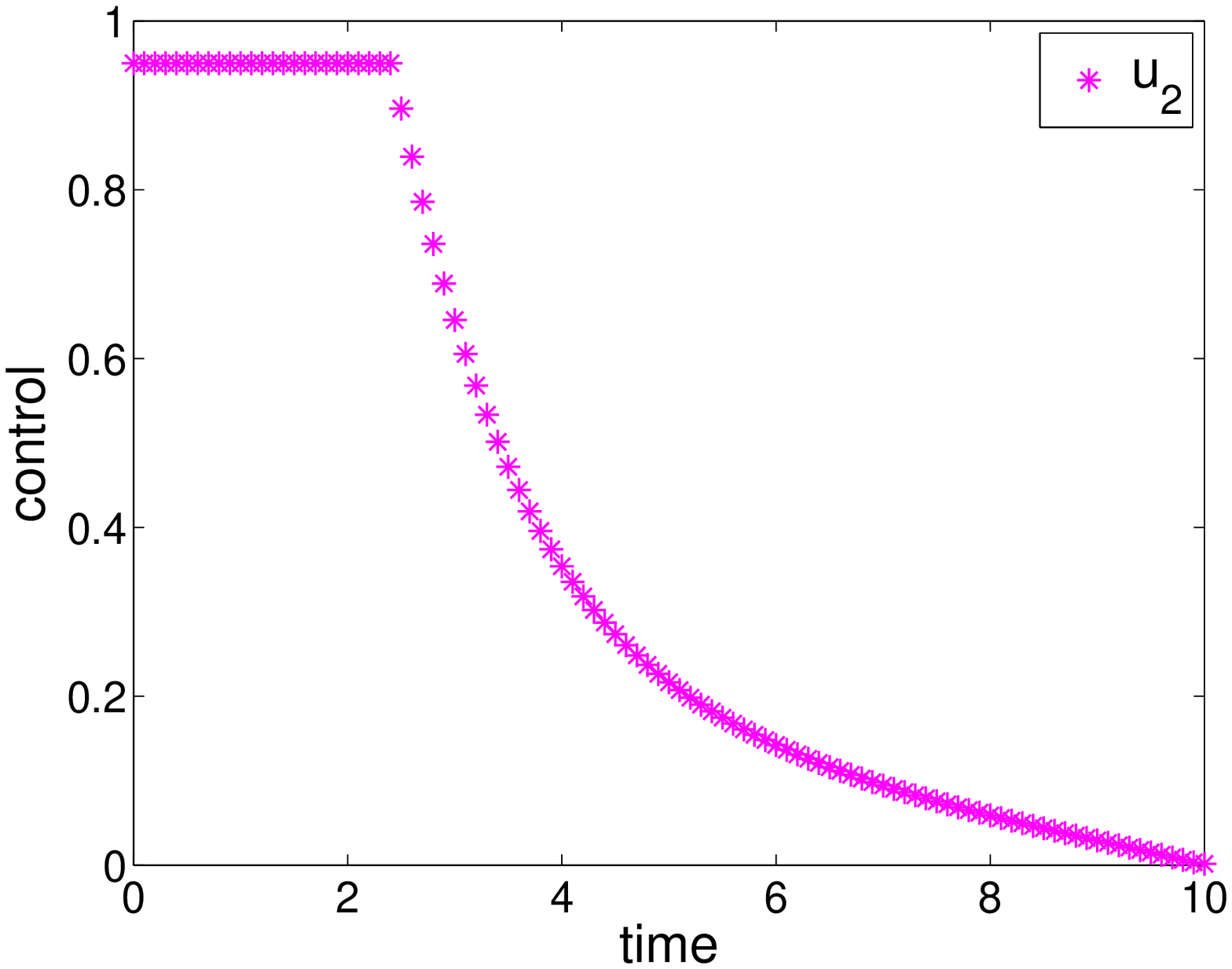,width = 0.33\textwidth}\label{fig:f3:control:MOP3}}
\subfigure[MOP1]{\epsfig{file=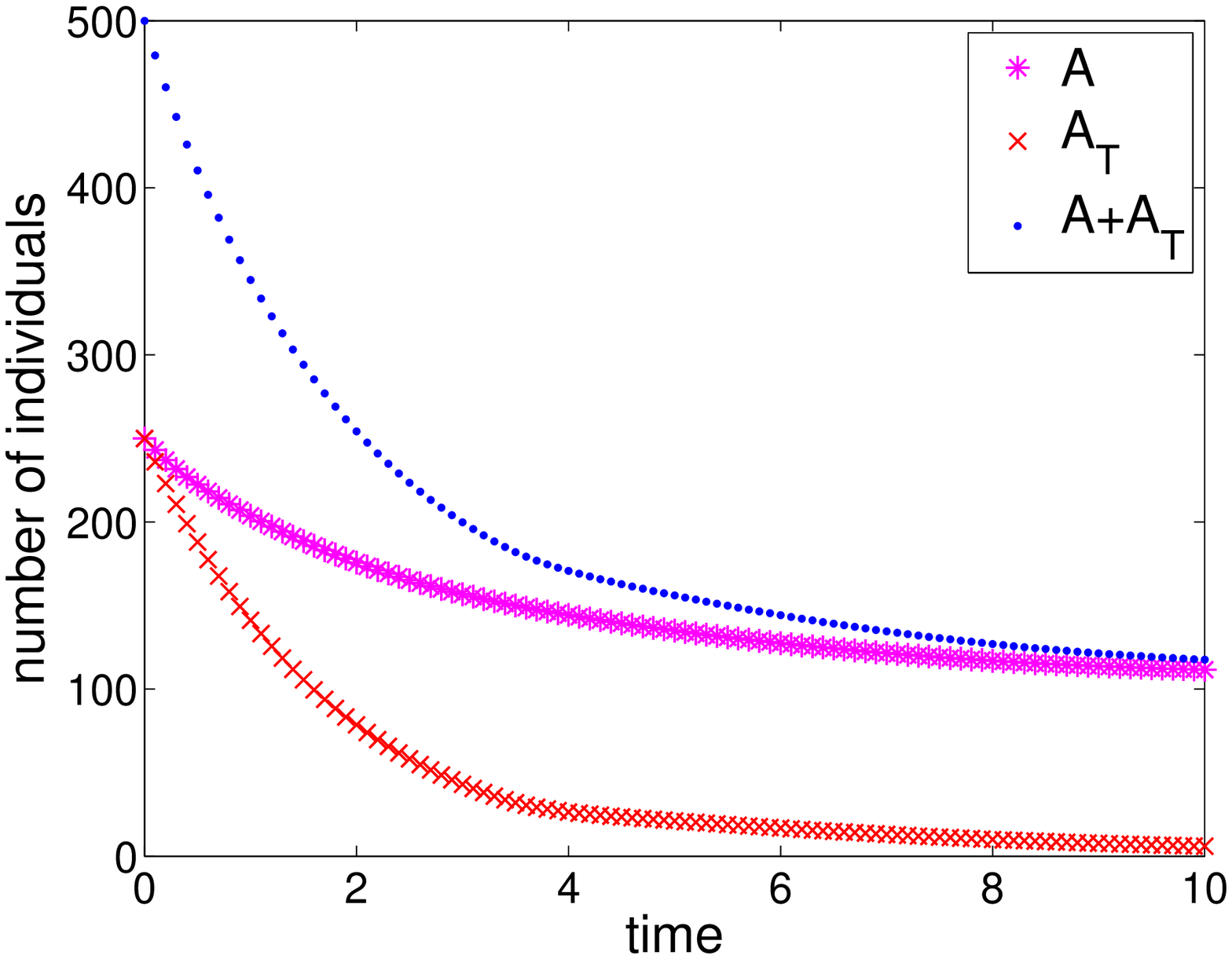,width = 0.33\textwidth}\label{fig:f3:inds:MOP1}}%
\subfigure[MOP2]{\epsfig{file=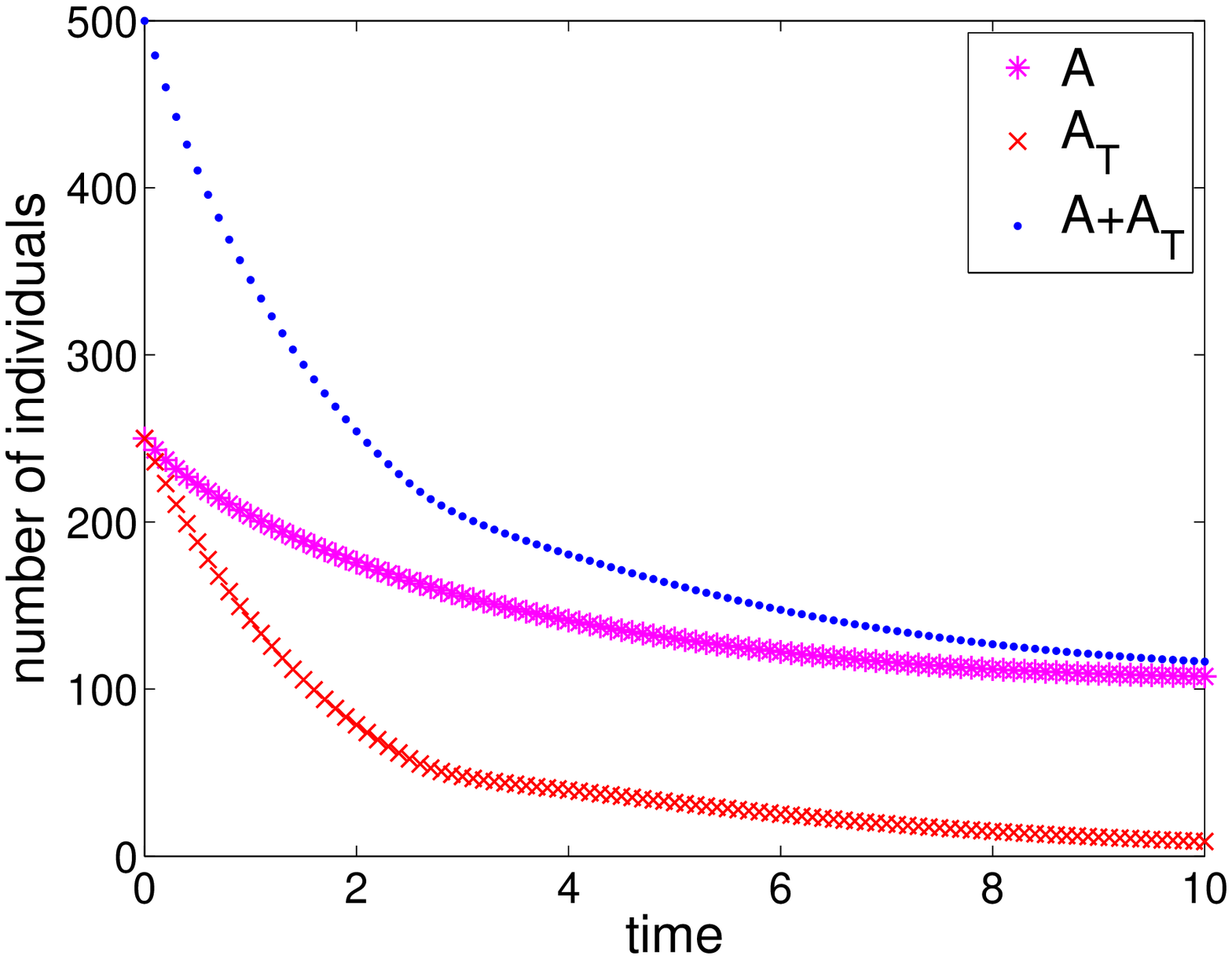,width = 0.33\textwidth}\label{fig:f3:inds:MOP2}}%
\subfigure[MOP3]{\epsfig{file=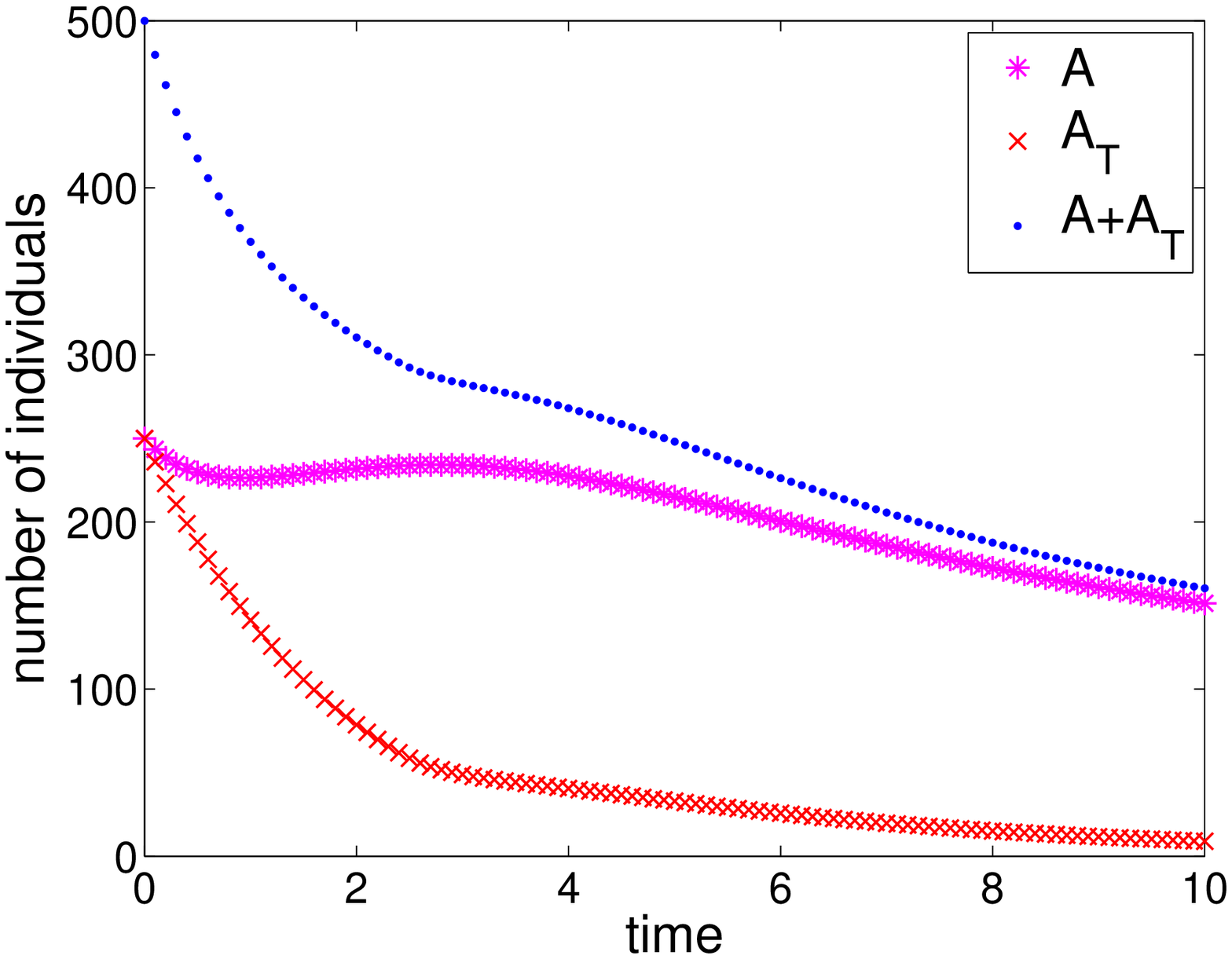,width = 0.33\textwidth}\label{fig:f3:inds:MOP3}}
\caption{Controls and infected individuals for $f_2 = 3$.}\label{fig:f3}
\centering
\subfigure[MOP1]{\epsfig{file=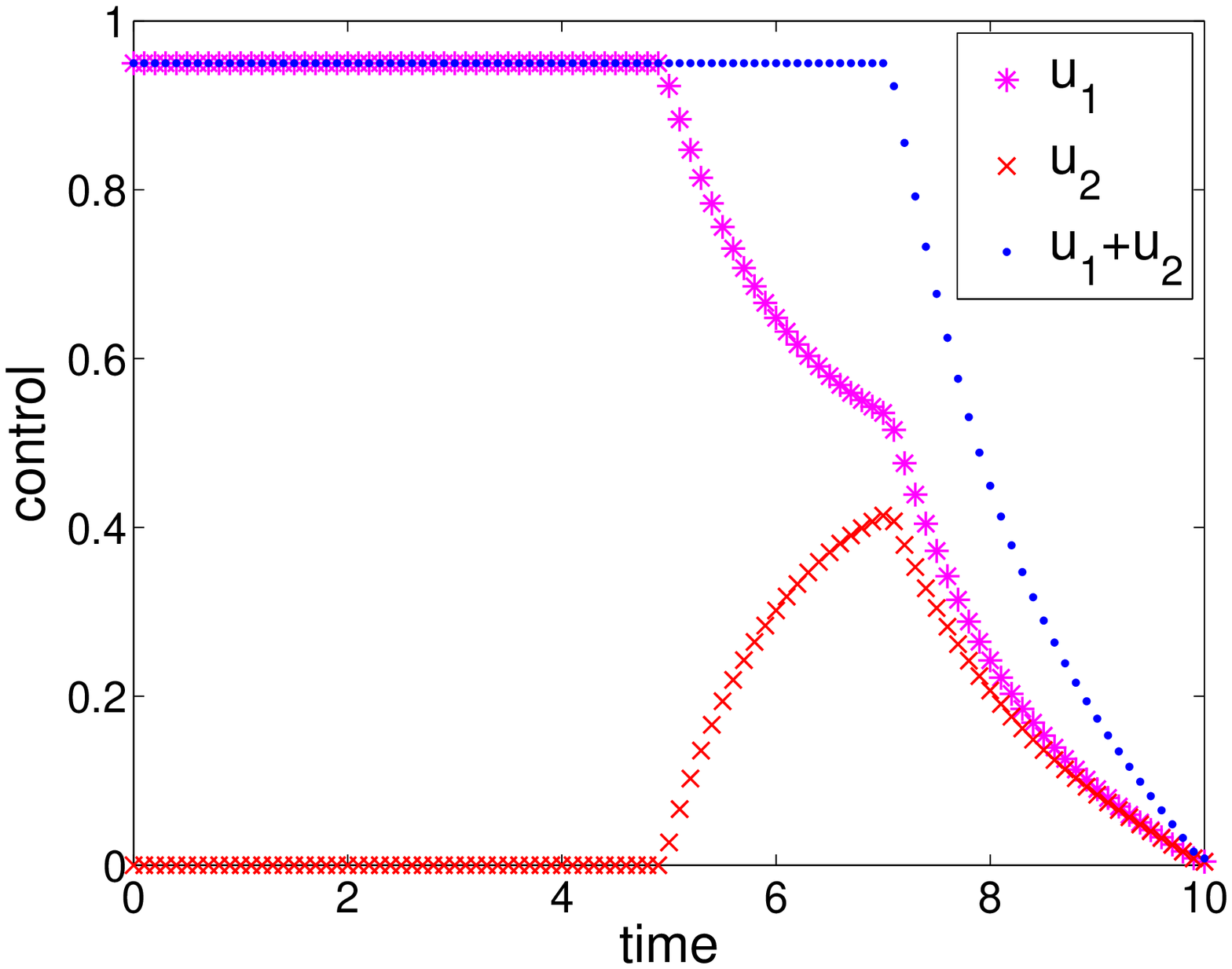,width = 0.33\textwidth}\label{fig:f6:control:MOP1}}%
\subfigure[MOP2]{\epsfig{file=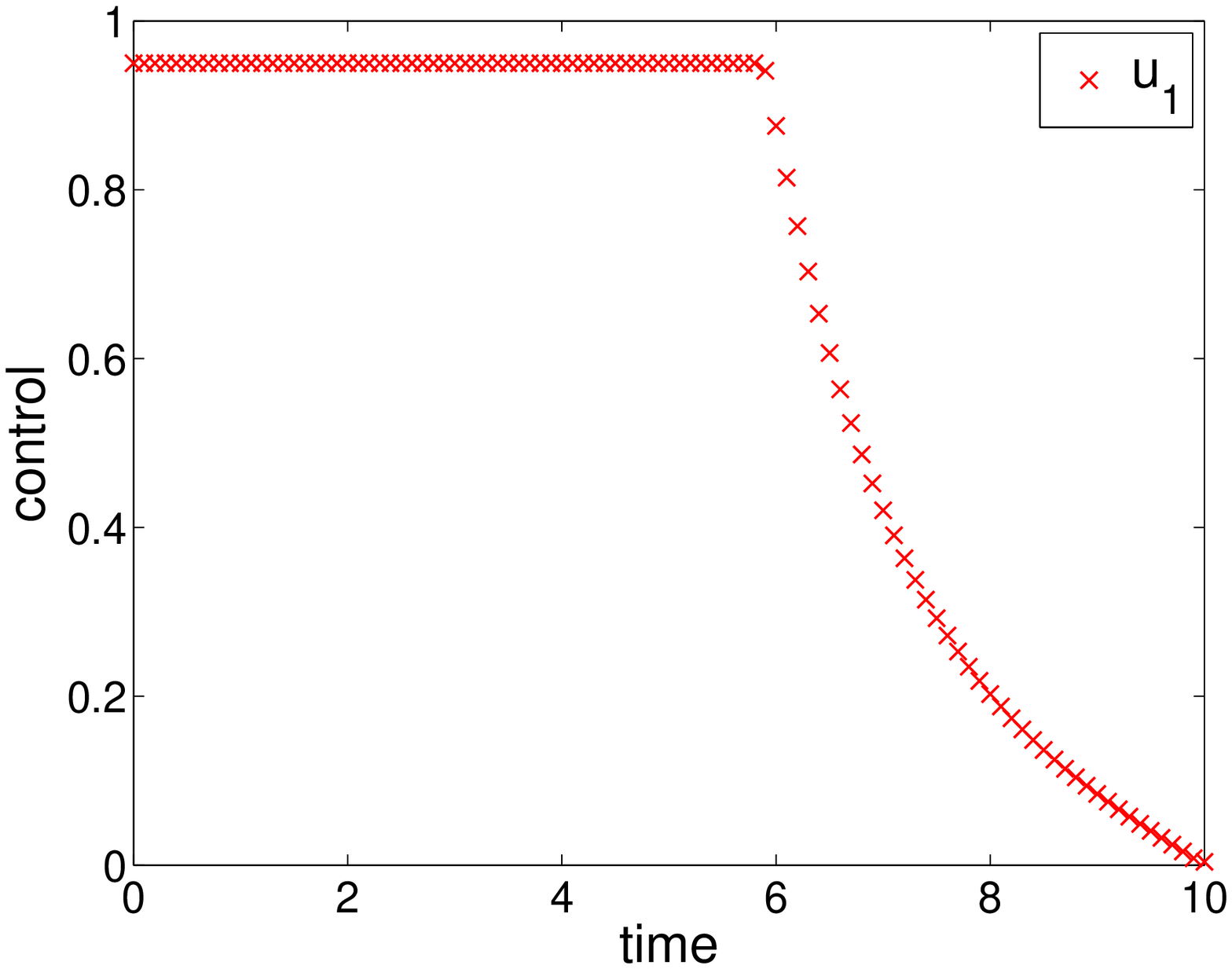,width = 0.33\textwidth}\label{fig:f6:control:MOP2}}%
\subfigure[MOP3]{\epsfig{file=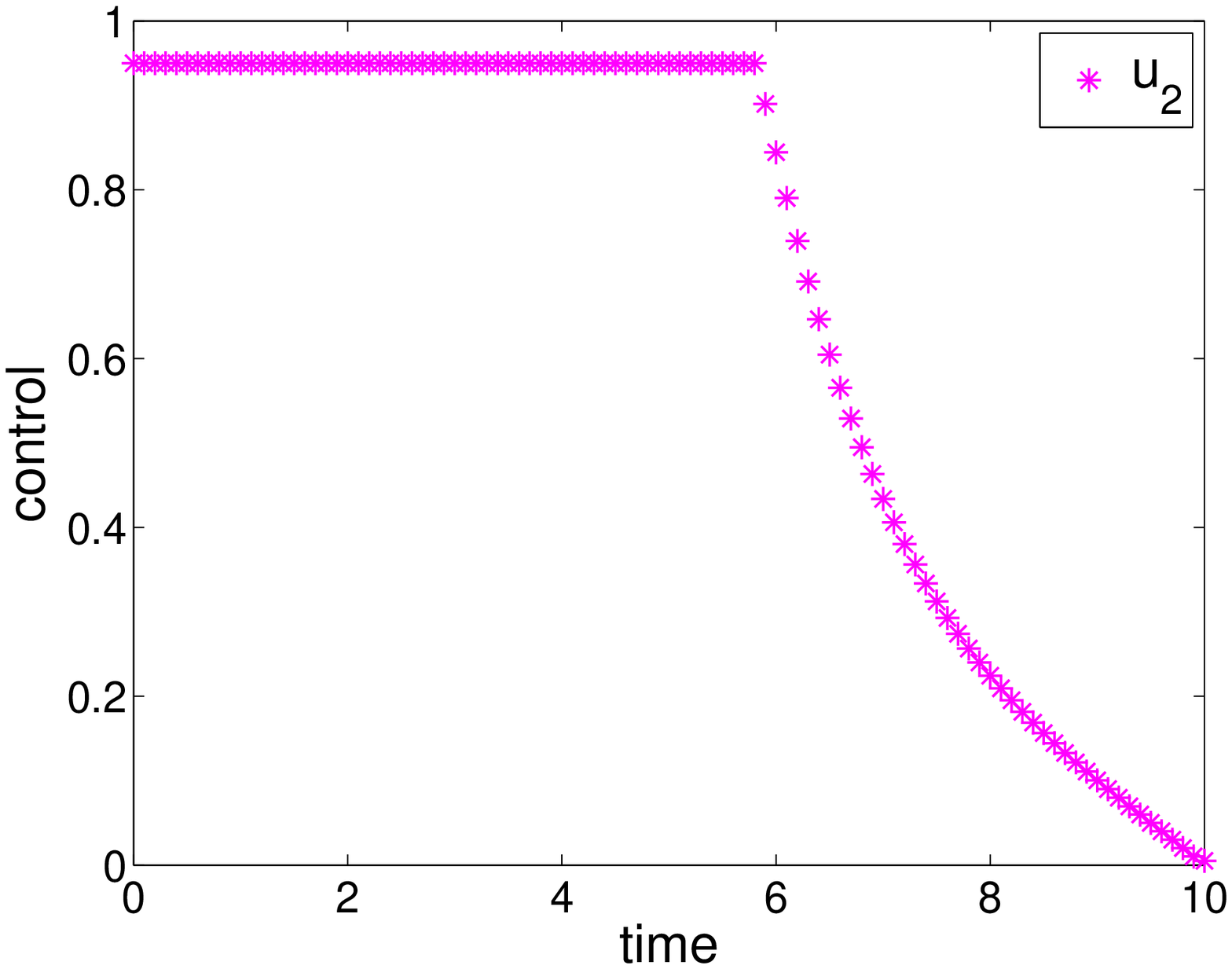,width = 0.33\textwidth}\label{fig:f6:control:MOP3}}
\subfigure[MOP1]{\epsfig{file=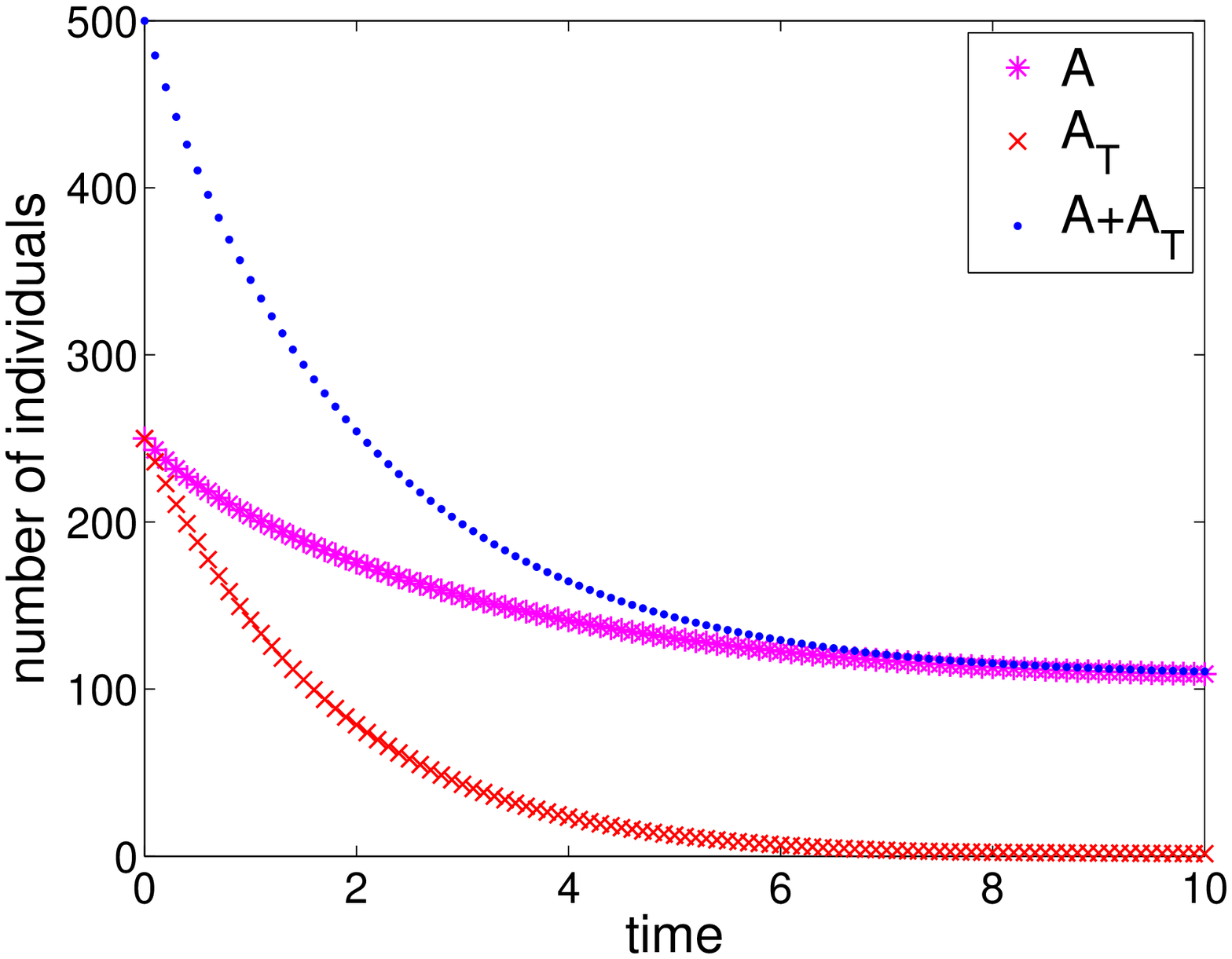,width = 0.33\textwidth}\label{fig:f6:inds:MOP1}}%
\subfigure[MOP2]{\epsfig{file=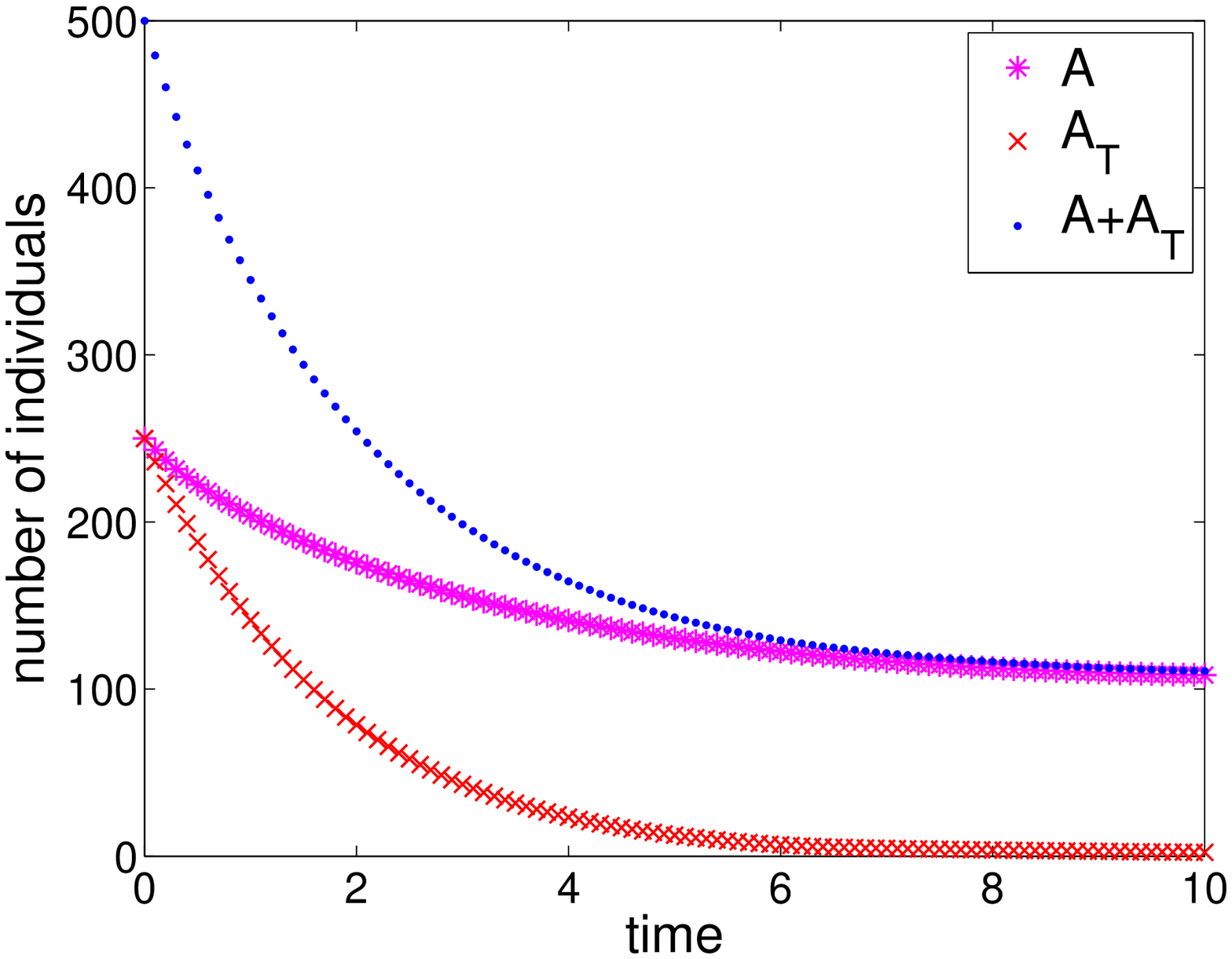,width = 0.33\textwidth}\label{fig::f6:inds:MOP2}}%
\subfigure[MOP3]{\epsfig{file=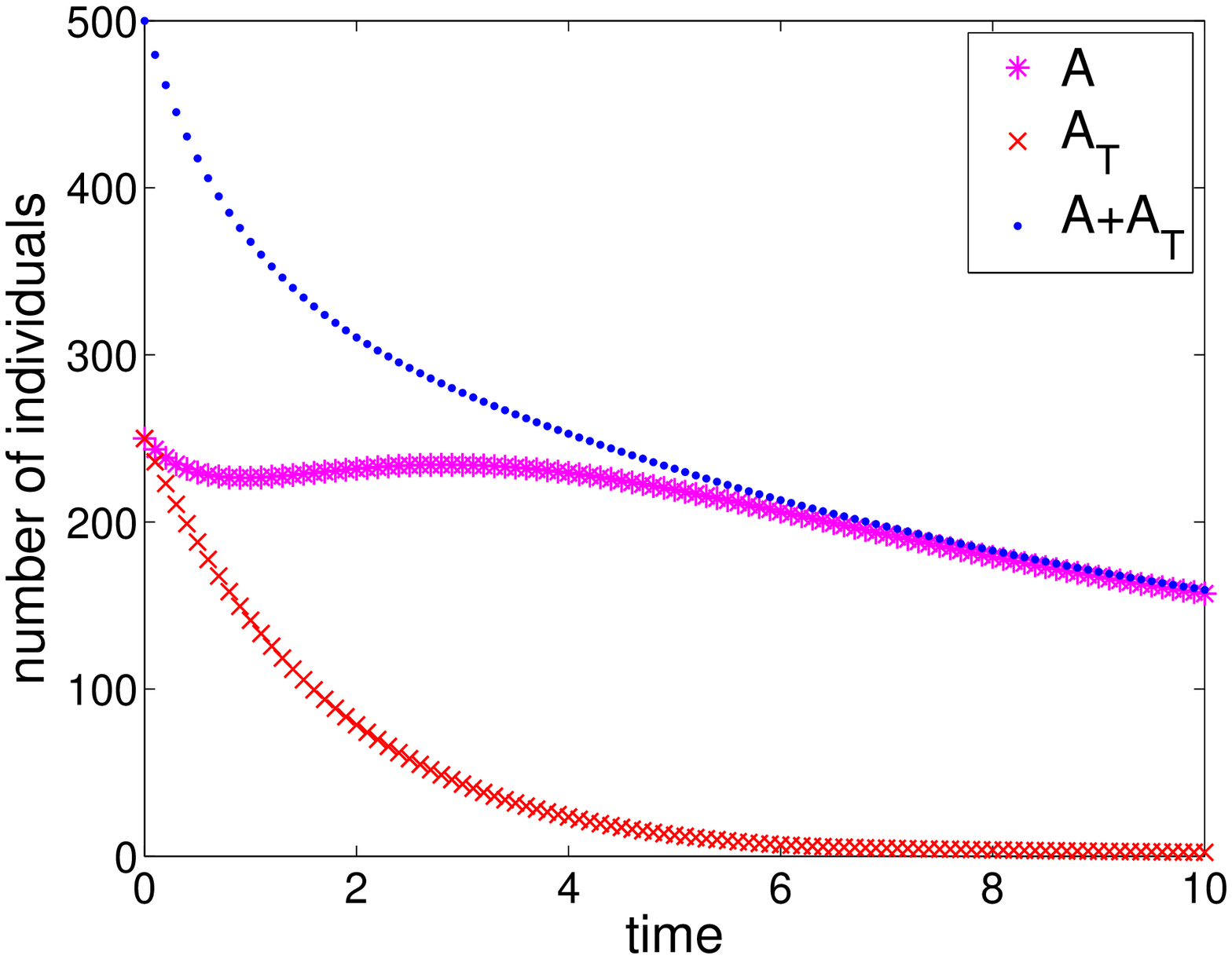,width = 0.33\textwidth}\label{fig:f6:inds:MOP3}}
\caption{Controls and individuals with AIDS and active TB for $f_2 = 6$.}\label{fig:f6}
\end{figure}
For the extreme solutions, i.e., those corresponding to the maximum and minimum
amounts of controls applied though the period of study, the dynamics of $A$,
$A_T$ and $(A+A_T)$ can be observed in Figures~\ref{fig:f0:inds}
and~\ref{fig:max:control}. Without applying the controls, as discussed above,
the corresponding dynamics are identical for MOP1, MOP2 and MOP3 and are
presented in Figures~\ref{fig:f0:inds}. This scenario corresponds to the natural
progression of the diseases. When only the control $u_2$ is applied, the number
of $A$ and $A_T$, as well as their sum, are larger than for cases when $u_1$
and $u_1+u_2$ are implemented. The dynamics for $A$ and $A_T$ are identical
for MOP1 and MOP3, which are depicted in Figure~\ref{fig:max:control:MOP12}.

To analyze intermediate scenarios, two solutions lying on intermediate regions
of the trade-off curves are selected. This is done as follows. The objective
space is divided by two horizontal lines corresponding to $f_2=3$ and $f_2=6$
(dashed lines in Figure~\ref{fig:pf}). The intersections of these lines with
the trade-off curves give the three different solutions, with each of them
corresponding to MOP1, MOP2 and MOP3. These solutions are selected for discussing
the dynamics of the diseases. These lines can be interpreted as the constraints
defining available resources for treatment. In turn, the selected solutions represent
the best treatment options in such circumstances, as they allow to achieve
the lowest values of $f_1$. It is worth noting that the solutions on different
curves, shown in Figure~\ref{fig:pf}, are identically distributed with respect
to the values of $f_2$ due to the use of the $\epsilon$-constrain method.
Since MOP1--3 were solved for the same values of $\epsilon$, the corresponding
solutions can be used for a fair comparison.
Figure~\ref{fig:f3} presents the trajectories of the control variables and the
dynamics of $A$ and $A_T$ for solutions corresponding to $f_2=3$.
From Figures~\ref{fig:f3:control:MOP1}--\ref{fig:f3:control:MOP3}, it can be seen
that the changes of the total amount of control measures are similar.
For MOP1, the total control is composed of $u_1$ and $u_2$,
which change differently during the period of study. There is a peak in $u_2$,
taking place between the third and forth years, after which it decreases.
The dynamics of $A$ and $A_T$ have similar trajectories, as shown
in Figures~\ref{fig:f3:inds:MOP1}--\ref{fig:f3:inds:MOP3}. However,
there is a minor increase in $A$ for MOP3 in the beginning of the forth year.
Similar trends for the control and state variables are observed in
Figure~\ref{fig:f6}. Though, the peak in $u_2$ occurs later with a higher value
and the number of individuals with AIDS and active TB is lower
due to the large amount of medication.

Solutions obtained using multiobjective optimization can provide comprehensive
insights about the optimal strategies and the diseases dynamics resulting from
implementation of those strategies. Since visual representations can help
to better understand results and spot patterns that are not obvious at first,
in what follows each of the variables $u_1$, $u_2$, $A$ and $A_T$ is defined
as a function of time and the objective to which it is conflicting. By doing so,
it is possible to provide the visualization of the entire optimal set of each
variable. The set of optimal values defines a surface. Slicing a surface gives
the trajectory of the corresponding dynamic over the period of study.
\begin{figure}
\centering
\subfigure[$u_1$]{\epsfig{file=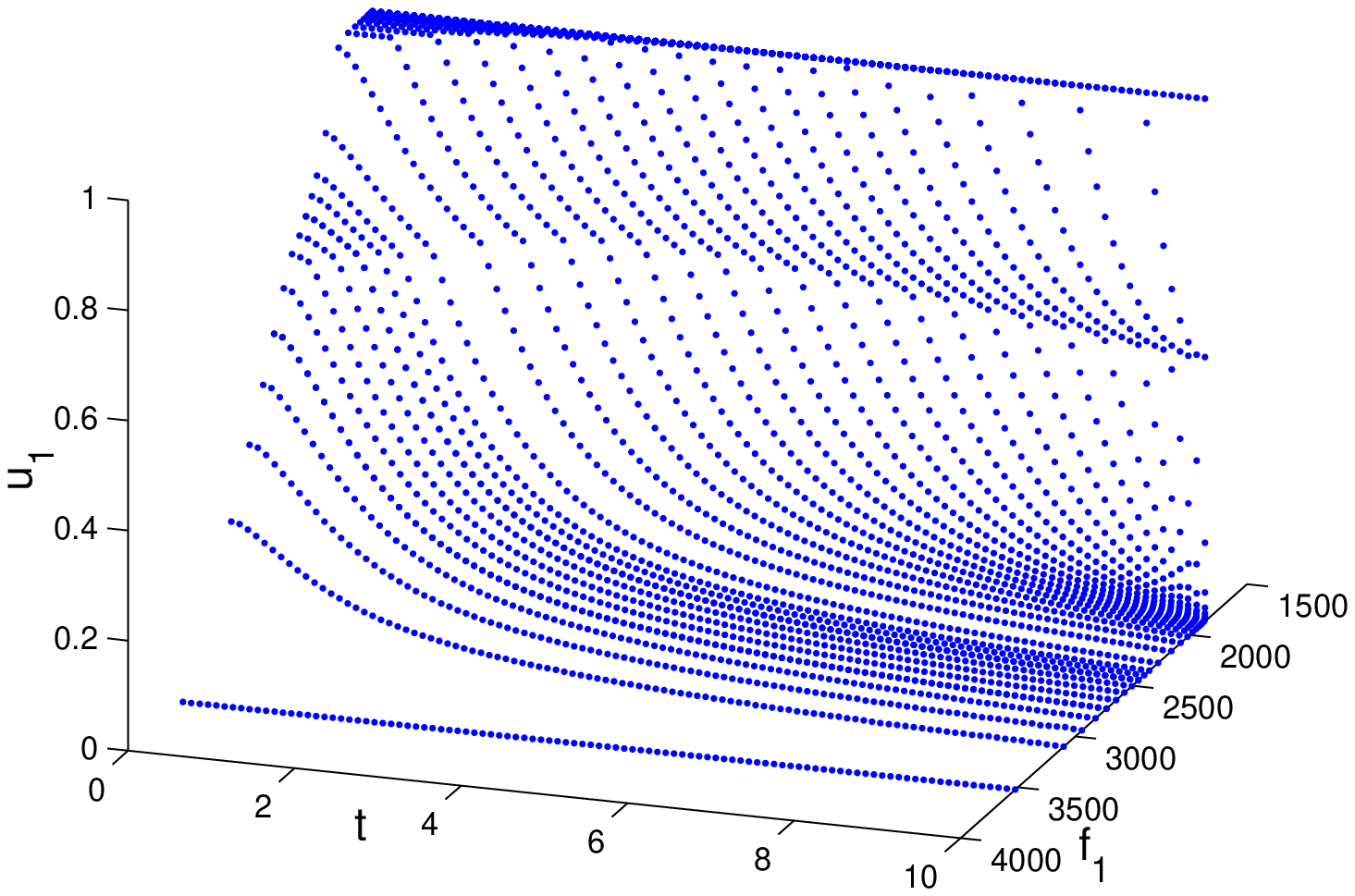,width = 0.45\textwidth}\label{fig:surface:controls:u1}}%
\subfigure[$u_2$]{\epsfig{file=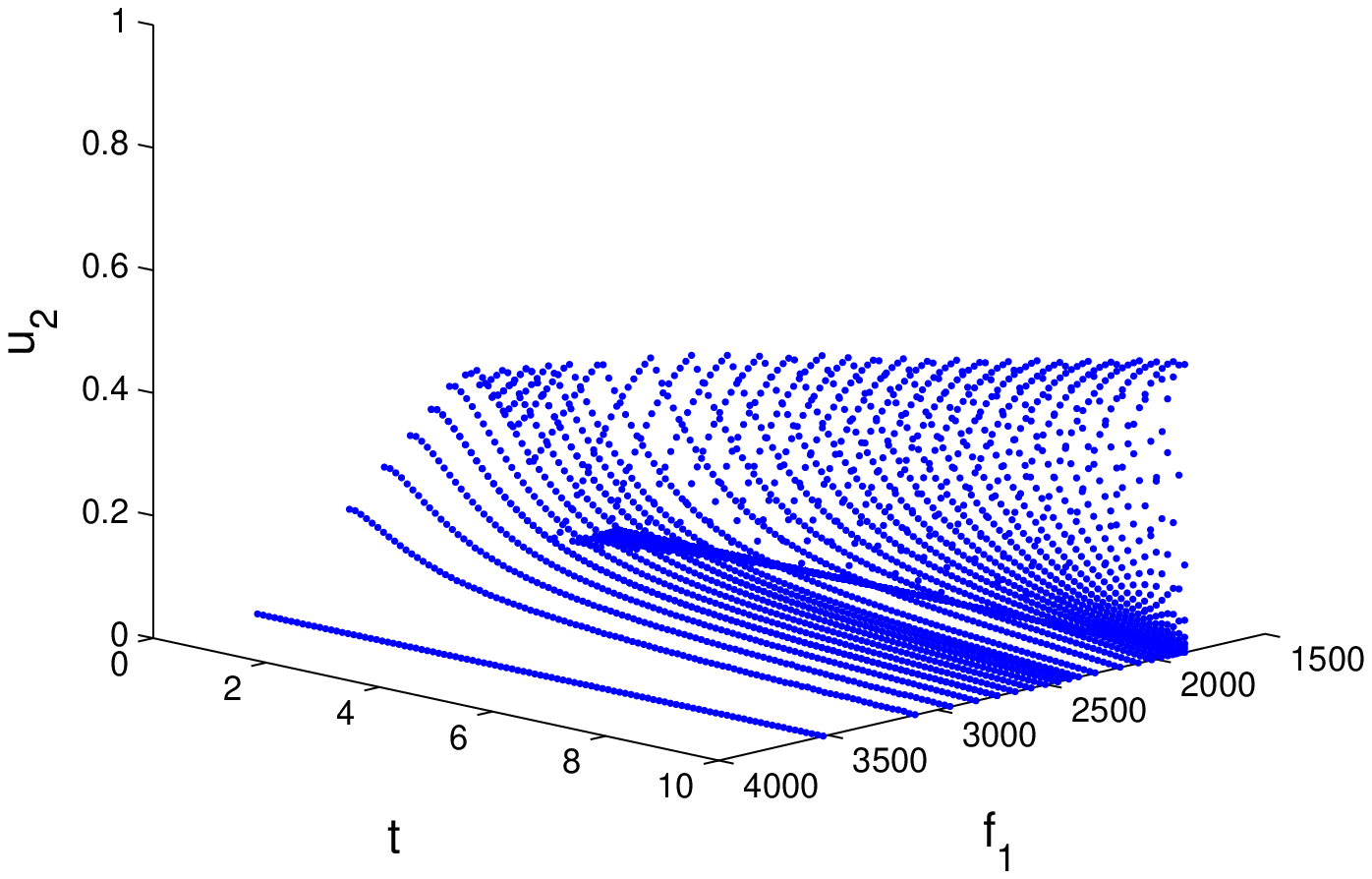,width = 0.45\textwidth}\label{fig:surface:controls:u2}}
\subfigure[$u_1+u_2$]{\epsfig{file=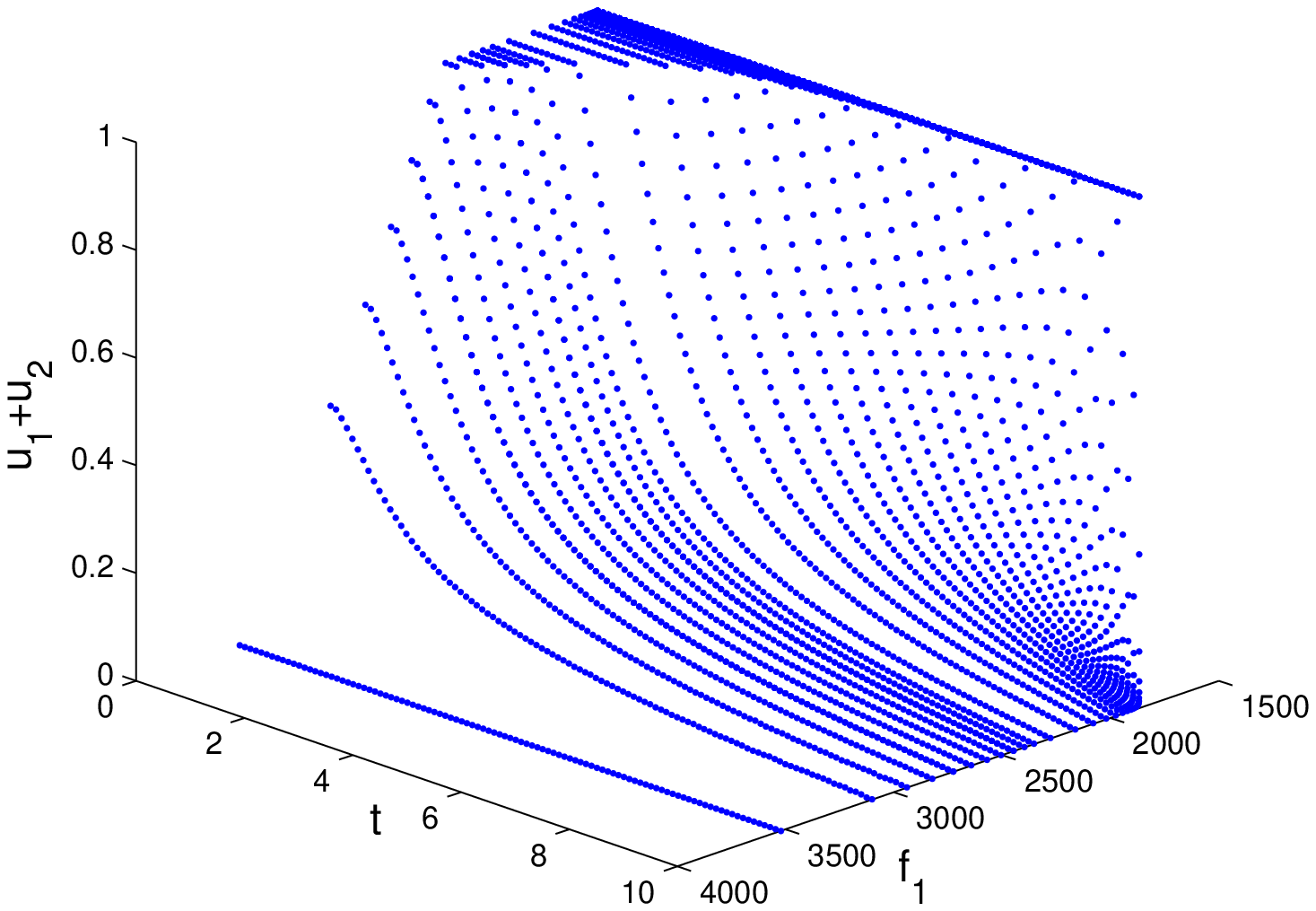,width = 0.45\textwidth}\label{fig:surface:controls:u1_u2}}
\caption{Discrete representation of the optimal surface for $u_1$, $u_2$ and $u_1 + u_2$.}\label{fig:surface:controls}
\centering
\subfigure[$A$]{\epsfig{file=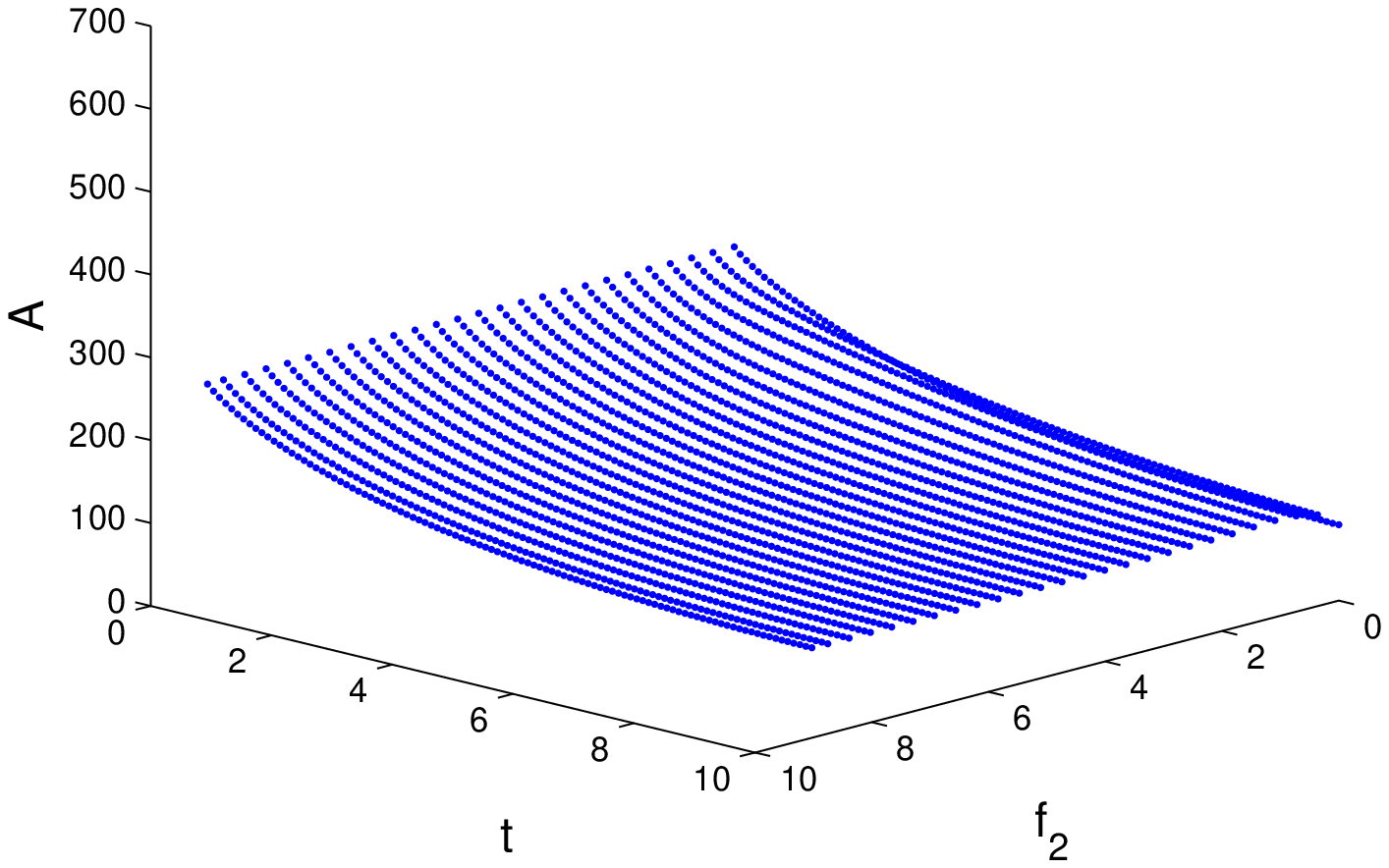,width = 0.45\textwidth}\label{fig:surface:inds:A}}%
\subfigure[$A_T$]{\epsfig{file=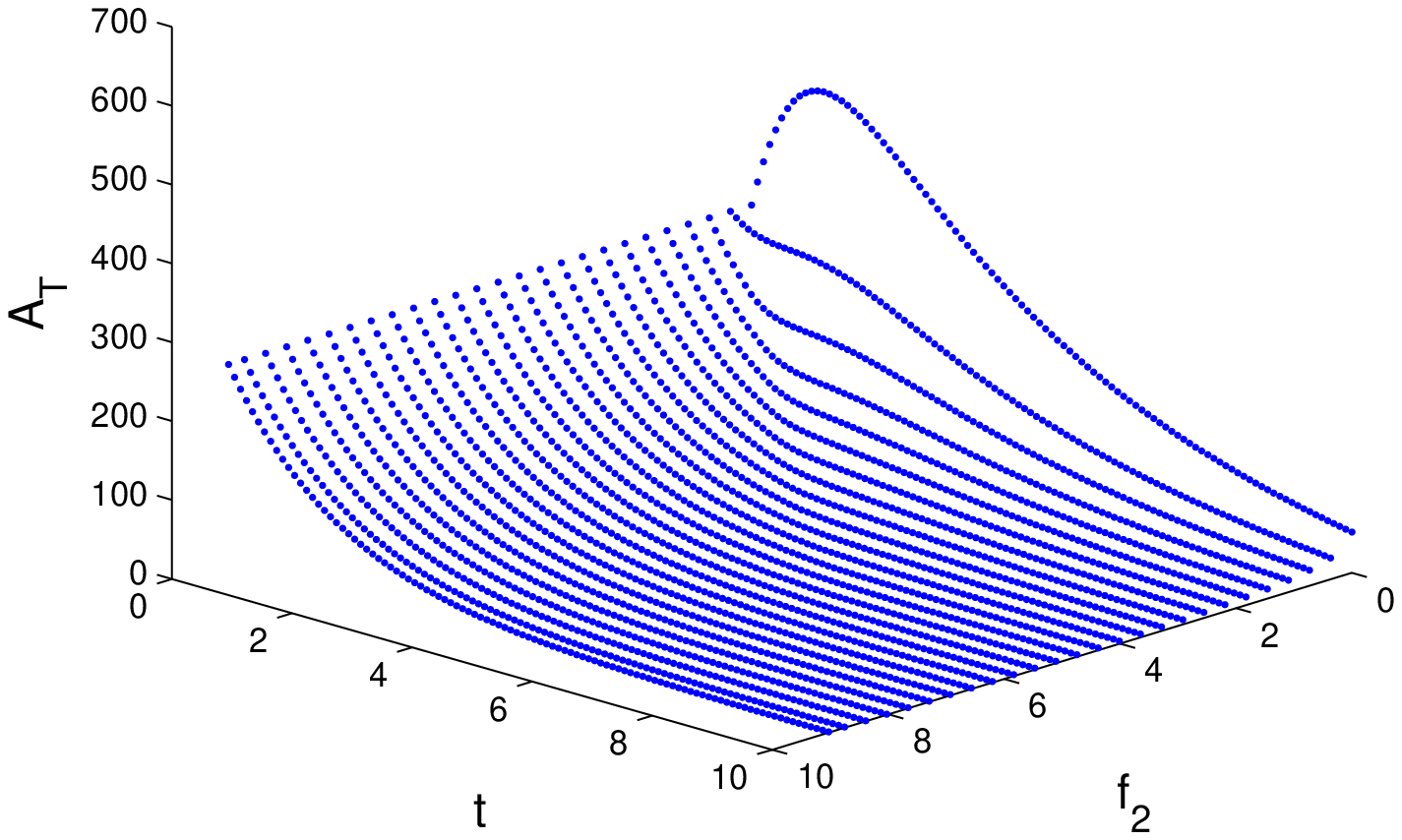,width = 0.45\textwidth}\label{fig:surface:inds:AT}}
\subfigure[$A+A_T$]{\epsfig{file=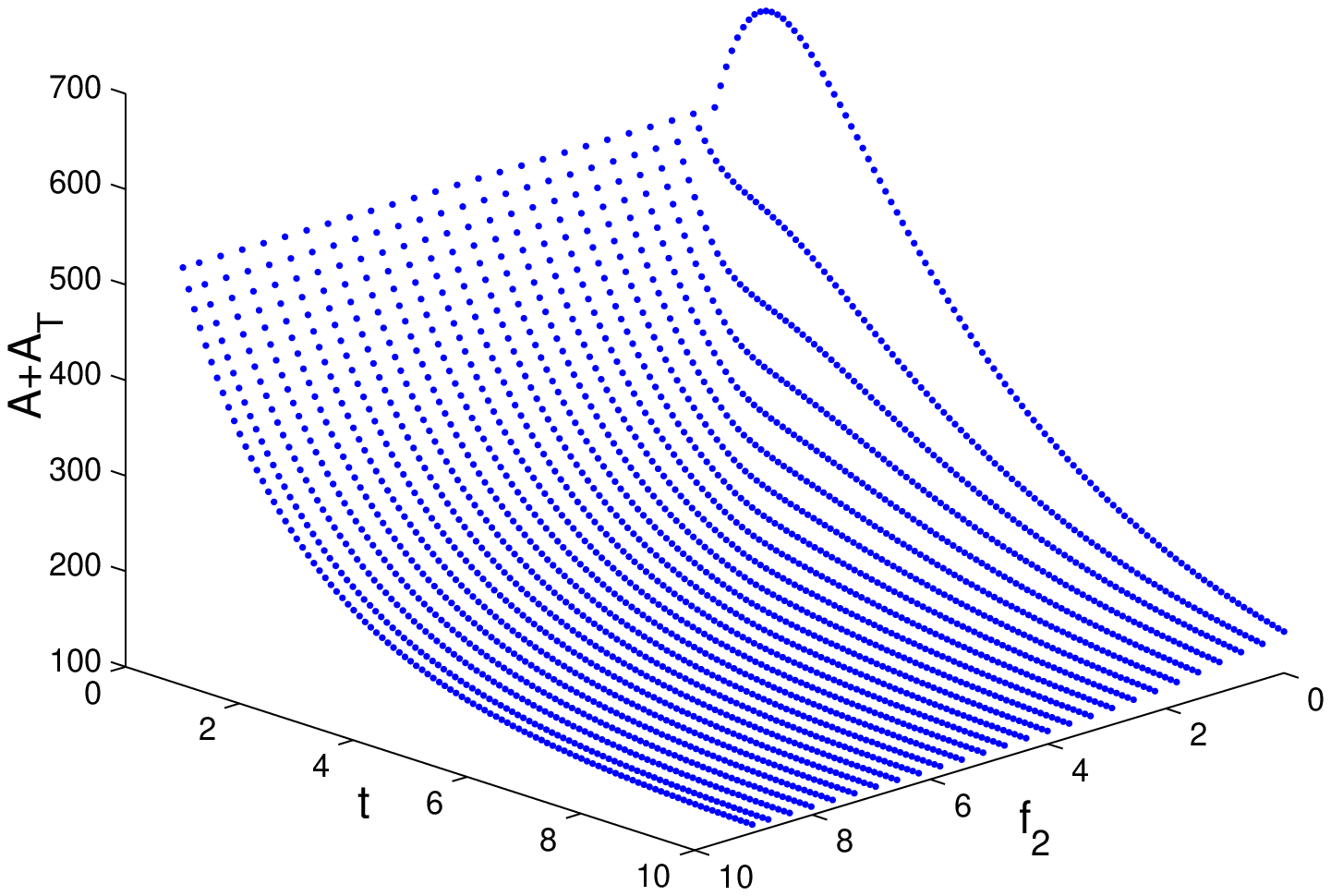,width = 0.45\textwidth}\label{fig:surface:inds:A_AT}}
\caption{Discrete representation of the optimal surface for $A$, $A_T$ and $A+A_T$.}\label{fig:surface:inds}
\end{figure}
Figure~\ref{fig:surface:controls} shows the discrete representations of surfaces
defined by the controls over the whole Pareto optimal region. On the other hand,
the discrete representations of surfaces determined by the responses of $A$
and $A_T$ to optimal treatment strategies are illustrated in
Figure~\ref{fig:surface:inds}. The plots for other classes
of human population can be obtained in a similar way.


\section{Conclusions}
\label{sec:conc}

This paper investigates a mathematical model for TB-HIV/AIDS coinfection
recently proposed in \cite{SiTo15}. A multiobjective formulation is proposed. 
This approach avoids the use of weight parameters and allows
to obtain a wide range of optimal control strategies, which offer useful
information for effective decision making. Two clearly conflicting
objectives are defined for search the optimal controls. The first objective
reflects aspirations in controlling TB and HIV/AIDS diseases, whereas the second
objective aims to reduce the costs of implementing control policies.
The present study extends the previous work~\cite{SiTo15} by the extensive
analysis of the optimal control in the TB-HIV/AIDS coinfection model,
which enriches the knowledge about the model. Indeed, it is important
not just to formulate a model but also to obtain useful information about
the process modeled. Our simulation results reveal the optimal treatment
strategies for TB and HIV infections and exposure to medication of different
fractions of the population. This can be used as an input for planning activities
to fight against TB and HIV. The choice of a final solution can be made including
the goals of public healthcare and available funds.
The results here obtained clearly demonstrate the usefulness
and advantages of a multiobjective approach. The presence of the clearly
conflicting objectives gives rise to a set of optimal solutions representing
different trade-offs between them. Thus, each obtained solution reveals
different perspectives on coping with AIDS and active TB diseases. The treatment
of individuals infected by both HIV and TB can provide the best effects, except
for the extreme scenarios. As analyses showed, the set of optimal trade-off
solutions can offer to the decision maker an understanding of all possible
trends in applying the controls. Moreover, the dynamics of different classes
of individuals in the population appears as a response to the implemented
treatment measures. The ability to obtain, analyze and choose from a set of
alternatives, constitutes the major advantage of the proposed approach,
motivating its practical use in the process of planning intervention measures
by health authorities. As future work, we intend to study the effects
of the parameters in the TB-HIV/AIDS coinfection model. Also,
it would be interesting to investigate the population dynamics resulting
from implementation of optimal treatment policies found by optimizing
various types of objectives. Considering $L^1$ objectives in the optimal
control problem is also the subject of future work.


\section*{Acknowledgements}

Silva and Torres were supported by Portuguese funds through 
the Center for Research and Development in Mathematics 
and Applications (CIDMA) and the Portuguese Foundation 
for Science and Technology (FCT),
within project UID/MAT/04106/2013; 
and by the FCT project TOCCATA, 
ref. PTDC/EEI-AUT/2933/2014.
Silva is also grateful to the FCT 
post-doc fellowship SFRH/ BPD/72061/2010.
The authors would like to thank two anonymous referees 
for valuable comments and suggestions.



\end{document}